\documentclass[a4paper,12pt,oneside]{amsart}
\usepackage[top=1in,bottom=1.1in,left=1.35in,right=0.95in]{geometry}
\usepackage{tikz}
\usepackage{amssymb}
\usepackage[utf8]{inputenc}
\usepackage[T1]{fontenc}
\usepackage[french,english]{babel}
\usepackage[all]{xy}
\usepackage{soul}
\usepackage{ulem}
\usepackage{shuffle}
\usepackage{hyperref}
\usepackage{stmaryrd}
\usepackage{graphicx}
\usepackage{array}
\usepackage{color}
\usepackage{xcolor}
\usepackage{slashed}
\usepackage{amsthm}
\usepackage{amsmath}
\usepackage{amssymb}
\usepackage{mathrsfs}
\usepackage{setspace}
\newtheorem{Lemme}{Lemma}[section]
\newtheorem{Prop}{Proposition}[section]
\newtheorem{Def}{Definition}[section]
\newtheorem{Rem}{Remark}[section]
\newtheorem{Thm}{Theorem}
\newtheorem*{Thmintro}{Theorem}
\newtheorem{Cor}{Corollary}[section]

\newtheorem{Exa}{Example}
\newtheorem{CExa}{Counter-Example}

\newcommand{\Preuve}{\noindent\textbf{Proof:\ }}
\newcommand{\eb}{\vspace{-1.5em}\begin{flushright}$\Box$\end{flushright}}
\newcommand{\eq}[1][r]
   {\ar@<-3pt>@{-}[#1]
    \ar@<-1pt>@{}[#1]|<{}="gauche"
    \ar@<+0pt>@{}[#1]|-{}="milieu"
    \ar@<+1pt>@{}[#1]|>{}="droite"
    \ar@/^2pt/@{-}"gauche";"milieu"
    \ar@/_2pt/@{-}"milieu";"droite"}

\newcommand{\cF}{{\mathcal{F}}}

\newcommand{\cT}{{\mathcal{T}}}

\newcommand{\Hol}{{\mathit{Hol}}}

\newcommand{\cH}{{\mathcal{H}}}
\newcommand{\bH}{\mathrm{\bf H}}
\newcommand{\bG}{\mathcal{G}}
\newcommand{\bT}{\mathrm{\bf T}}

\begin{document}
\selectlanguage{english}
\title{On the pseudogroup of local transformations commuting with a transversely elliptic operator and the existence of transverse metric}
\author{Wenran Liu}
\address{Wenran Liu, Sino-European Institute of Aviation Engineering, Civil Aviation University of China, 2898 Jinbei Street Dongli District, 300300 Tianjin, P.R.China }
\email{wrliu@cauc.edu.cn}
\thanks{This research project is supported by the Starting Research Fund from Civil Aviation University of China, Project No. 2020KYQD58}
\date{}

\begin{abstract}
The group of diffeomorphisms commuting with an elliptic operator on a manifold is a compact Lie group under the Compact-Open topology. In foliation theory,
pseudogroups were introduced by Sacksteder. The pseudogroup of local transformations commuting with a basic differential operator possesses the equicontinuity and the quasi-analyticity properties when conditions on the operator are given. These properties serve to construct a transverse metric on the normal bundle under a good condition on the operator. For this, the Average Method is applied as in the construction of basic connections on foliated bundles.
\end{abstract}

\maketitle

\tableofcontents
\normalsize

\begin{spacing}{1.1}
\section{Introduction}
The notion of bundle-like metric on foliated manifold was first introduced by Reinhart \cite{Rei}. A foliated manifold equipped with a bundle-like metric is called a Riemannian foliation. An equivalent definition is the existence of a transverse (holonomy-invariant) metric on the normal bundle of the foliation. In \cite{Mol0, Mol}, Molino proposed a theory for the study of Riemannian foliation. In \cite{Car} and in \cite{Bre}, Carrière and Breuillard-Gelander studied the growth of Riemannian foliations.  In \cite{Ka}, El.Kacimi-Alaoui introduced the notion of transversely elliptic operator and proved the association of a transversely elliptic operator over a Riemannian foliation and an equivariant operator over a compact manifold constructed by Molino's theory. In \cite{BKR}, J.Bruning, F.Kamber and K.Richardson proposed an index theory for basic Dirac operators on Riemannian foliations. In \cite{GL}, A.Gorokhovsky and J.Lott gave an index formula for a transverse Dirac-type operator in case of Abelian Molino sheaf. As their work are all based on the Riemannian assumption of the foliation, El.Kacimi-Alaoui asked in \cite[Question2.9.2]{AAHR} whether a transversely elliptic operator exists only if the foliation is Riemannian.\smallskip\\
\indent In the present work, let us answer this question by proposing a sufficient condition on $P$. Let us reformulate this question as follows:
\smallskip
\begin{center}
\begin{minipage}{0.8\textwidth}
Let $M$ be a compact manifold with a foliation $\cF$. Suppose that there exists a basic transversely elliptic operator $P\!:\!C^\infty(M)\!\rightarrow\! C^\infty(M)$ of order $m\!\geq\! 2$ acting on the complex-valued smooth functions $C^\infty(M)$ of $M$. Is the foliation $\cF$ Riemannian?
\end{minipage}
\end{center}
\smallskip

Here the word ``basic'' means that at each point of $M$, there exists a local foliated chart $U\!\simeq\!(x,y)$ $($locally $\cF$ being defined by $dy\!=\!0$$)$ such that over $U$ the operator $P$ has the following expression $P\!=\!\sum_{|s|=0}^m a_s(y)\frac{\partial^{|s|}}{\partial y_1^{s_1}\cdots\partial_n^{s_n}}$ where $s\!=\!(s_1,\cdots,s_q)\in\mathbb{N}^q$, $|s|\!\!=\!\!s_1\!+\!\cdots\!+\!s_q$ ($q$ being the codimension of $\cF$) and $a_s$ are smooth functions of $y$.\smallskip\\
\indent In \cite{Liu}, the \textbf{Average Method} is applied for the construction of a basic connection on a foliated bundle over a Riemannian foliation. Let us take a brief review. Let $\cT$ be a complete transversal of $\cF$ and consider the holonomy groupoid $\mathrm{Hol}(M,\cF)_\cT^\cT$ restricted to $\cT$. The closure $\overline{\mathrm{Hol}(M,\cF)_\cT^\cT}$ is well-defined to be a Lie, Hausdorff and proper groupoid, equipped with a Haar system. Together with ``cut-off'' functions, a basic connection is constructed by averaging a foliated one. \smallskip\\
\indent We hope that this method can also be applied to the construction of a transverse metric on $\nu\cF$, \emph{i.e.} a holonomy-invariant Riemannian metric on $\cT$. Unfortunately, since we do not assume $\cF$ to be Riemannian in our setting, $\overline{\mathrm{Hol}(M,\cF)_\cT^\cT}$ is not well-defined. Alternatively, we consider naturally the \textbf{holonomy pseudogroup} \textit{Hol} on $\cT$ introduced by Haefliger \cite{Hae1}.\smallskip\\
\indent In case that $\cF$ is Riemannian, in \cite{AC1}, \'{A}lvarez and Candel discussed the generalization from Riemannian foliations to equicontinuous foliated spaces, see also Tarquini's paper \cite{Tar}. In \cite[Appendix E]{Mol} (also in \cite{Kel}), Ghys pointed out that equicontinous foliated spaces can be considered as the ``topological Riemannian foliations'' and therefore the  In \cite{Hae,Hae2}, Haefliger discussed the quasi-analytic properties of \textit{Hol} and the existence of the closure of \textit{Hol}. More than the holonomy pseudogroup, in the present work, we consider the \textbf{pseudogroup} $\bG(P)$ of \textbf{local transformations} on $\cT$ \textbf{commuting with} $P$ and study the properties of equicontinuity and quasi-analyticity of $\bG(P)$.\smallskip\\
\indent In \cite{Fu}, Furutani studied the transformation group of a compact manifold commuting with an elliptic differential operator. He also did further works in \cite{Fu1} for an elliptic pseudo-differential operator. The main result of \cite{Fu} is that equipped with the Compact-Open topology, the group of diffeomorphisms commuting with an elliptic operator on a compact manifold is a compact Lie group using global analysis on manifolds. Although a global condition of $P$ such as the transverse ellipticity implicates the equicontinuity of $\bG(P)$, it does not allow to obtain sufficient local properties for a further construction. Clearly, some local conditions on $P$ are necessary. We will see in \textbf{Example} 2 that the present work is a \textbf{generalization} of \cite{Fu}. Indeed, when the foliated manifold is constructed by suspension, local conditions are unnecessary and the transverse ellipticity of $P$ is sufficient to obtain the existence of a transverse metric on $\nu\cF$.\smallskip\\
\indent In a local chart, $P$ can be regarded as a differential operator on an open subset $\Omega\!\subset\!\mathbb{R}^{\dim M}$. It is natural to suppose that $P$ has analytic coefficients, even more constant coefficients. Fortunately, these conditions implicate properties of $\bG(P)$ as the following tabular shows:
\begin{center}
\begin{tabular}{|c|c|}
\hline
Condition of $P$&Property of $\bG(P)$\\
\hline
Transversely elliptic&Equicontinuity\\
Analytic with triangular $1$-part&Quasi-analyticity\\ 
Constant coefficients&Transitivity\\ 
\hline
\end{tabular}
\end{center}
\indent These allow us to define a pseudogroup whose closure allows to apply the average method for the construction of transverse metric.   \\
The main result of this article is as follows:     
\begin{Thmintro}
Let $(M,\cF)$ be a compact foliated manifold and $P\!:\!C^\infty(M)\!\rightarrow\! C^\infty(M)$ be a basic differential operator, $\cT$ be a complete transversal satisfying \textrm{\bf Assumption} and $\bG(P)$ be the commuting pseudogroup.
\begin{enumerate}
\item If $P$ is transversely elliptic, $\bG(P)$ is equicontinuous with respect to the induced Riemannian distance;
\item If $P$ is transversely elliptic and analytic with triangular $1$-part at $z\!\in\!\cT$, $\bG(P)$ is quasi-analytic around $z$;
\item A sub-pseudogroup $\cH$ of $\bG(P)$ is quasi-analytic if its restriction to an open subset of $\cT$ meeting all $\cH$-orbits is quasi-analytic.
\end{enumerate}
\end{Thmintro}
In particular, when $P$ is basic transversely elliptic whose expression has constant coefficients in a local chart of $z$, the condition (2) of the Theorem above is satisfied. Moreover, when $P$ has constant coefficients on a complete transversal $\cT$ of $\cF$ satisfying \textbf{Assumption}, we have the following theorem.
\begin{Thmintro}
Let $(M,\cF)$ be a compact foliated manifold and $P\!:\!C^\infty(M)\!\rightarrow\! C^\infty(M)$ be a basic transversely elliptic operator. The foliation $\cF$ is a Riemannian foliation if there exists a complete transversal $\cT$ of $\cF$ satisfying \textbf{Assumption} where $P$ has constant coefficients.
\end{Thmintro}
\noindent In Section 2, we review some basic definitions and properties of pseudogroup, foliation, operator on foliated manifold. Under an assumption on transversal, we define the pseudogroup of local transformations commuting with $P$.\medskip\ \\
\noindent In Section 3, we impose conditions on $P$ such that $\bG(P)$ possesses sufficient properties for a further work.\medskip\ \\
\noindent In Section 4, we work under the condition that $P$ has constant coefficients. We first discuss the local/global transitivity of $\bG(P)$, secondly we define a sub-pseudogroup $\cH$ of $\bG(P)$ containing the holonomy pseudogroup and pass to its closure $\overline{\cH}$. Finally, a transverse metric is constructed with the help of $\overline{\cH}$.\\

\noindent {\bf\em Acknowledgements:} \\[0.5em]
\noindent The author is grateful to Paul-Emile Paradan, Moulay Tahar Benameur and Alexandre Baldare for various discussions on this project. This project is supported by the Starting Research Fund from Civil Aviation University of China, Project No.2020KYQD58.\\[0.5em]
\noindent The author welcomes and gratefully accepts all remarks, comments, doubts and criticisms regarding this paper. \\[1em]
\noindent Contact: \textbf{wrliu@cauc.edu.cn} also \textbf{wenran.liu@gmail.com}
\section{Prelimiaries}
\subsection{Pseudogroups of transformation}\ \\ \\
\indent In this section, we gather standard notions about pseudogroups, see \cite{Hae,Hae1,Sac,Sal1} for more details. All manifolds and all maps are assumed to be of class $C^\infty$. Let $\bT$ be a smooth manifold. Let us recall that a \textbf{local transformation} on $\bT$ is a diffeomorphism whose domain and range are connected open subsets of $\bT$. 
\begin{Def}[Haefliger \cite{Hae}, Salem\cite{Sal1}]\label{pseudogroup} A \textbf{pseudogroup} $\bH$ on $\bT$ is a collection of diffeomorphisms $h\!:\!U\!\rightarrow\!V$ $($$U,V$ being open subsets of $\bT$$)$ such that
\begin{itemize}
\item[$\bullet$] if $h\in \bH$, then its inverse $h^{-1}\in\bH$;
\item[$\bullet$] if $h:U\rightarrow V,h^\prime:U^\prime\rightarrow V^\prime$ belong to $\bH$, then the composition $h^\prime\! \circ\! h$ 
$$
U\cap h^{-1}(U^\prime)\stackrel{h}{\longrightarrow}V\cap U^\prime \stackrel{h^\prime}{\longrightarrow} h^\prime(V)\cap V^\prime;
$$
belongs to $\bH$;
\item[$\bullet$] the identity map of $\bT$ belongs to $\bH$;
\item[$\bullet$] if $h:U\rightarrow V$ belongs to $\bH$, then its restriction to each open subset of $U$ also belongs to $\bH$;
\item[$\bullet$] if $h:U\rightarrow V$ is a diffeomorphism between open subsets of $\bT$ which coincides in a neighborhood of each point of $U$ with an element of $\bH$, then $h\in\bH$.
\end{itemize}
\end{Def}
\indent A \textbf{sub-pseudogroup} of $\bH$ on $\bT$ is a pseudogroup on $\bT$ contained in $\bH$. Denote respectively $\operatorname{dom}h,\operatorname{im}h$ for the domain and range of $h\in \bH$. The set $\bH\!\cdot\! z\!=\!\big\{h(z)\mid h\!\in\! \bH,z\!\in\!\operatorname{dom}h\big\}$ is the \textbf{orbit} of a point $z\!\in\!\bT$. The \textbf{restriction} of $\bH$ to an open subset $\mathrm{\bf U}\subset \bT$ is the pseudogroup $\bH\vert_{\mathrm{\bf U}}\!=\!\big\{ h\in\bH\mid \operatorname{dom}h\cup\operatorname{im}h\!\subset \!\mathrm{\bf U} \big\}$.
The pseudogroup \textbf{generated} by a set $S$ of diffeomorphisms between open subsets of $\bT$ is the intersection of all pseudogroups that contains $S$. We shall say that a Riemannian metric $g$ on $\bT$ is $\bH$-\textbf{invariant} if for every $ h\in\bH$, every $z\in\operatorname{dom}h$, 
$$
h^*g\vert_{h(z)}\!=\!g\vert_z
$$ 
where the map $h^*$ is the standard map induced by $h$.      
\begin{Rem}
\'{A}lvarez Lopez and Moreira Galicia \cite{AG} focused on topological properties of pseudogroup. Thus in their work, $\bT$ is a topological space.
\end{Rem}
Let $\mathrm{Loct}(\bT)$ be the pseudogroup of all local transformations on $\bT$. It is equipped with the \textit{Compact-Open topology on Partial Maps with Open Domains}, see \cite{AB}. Every pseudogroup on $\bT$ is a sub-pseudogroup of $\mathrm{Loct}(\bT)$ equipped with the induced topology.
\subsection{Foliated manifold}\ \\ \\
\indent The needed basic concepts from foliation theory can be found in \cite{CC,Hec,Hec1}. Let $\cF$ be a foliation on a \textbf{compact} manifold $M$. We shall denote by $p\!=\!\dim\cF$ and $q\!=\!\operatorname{Codim}\cF$.
Let $T\cF$ denote the subbundle of vectors tangent to the leaves and $\nu\cF\!=\!TM\slash T\cF$ the normal bundle of $\cF$. 
\begin{Def} 
A \textbf{transversal} $\cT$ of $\cF$ is a submanifold of $M$ of dimension $q$ such that $\forall\, z\!\in\!\cT$, $T_z\cT\cap T\cF\vert_z\!=\!\{0\}$. A \textbf{complete transversal} is an immersed $($not necessarily connected$)$ submanifold, which intersects any leaf in at least one point.
\end{Def}
On $M$, there exists a complete transversal with a finite number of connected components. Indeed, take a cover of $M$ by foliated charts $\{U_i\}_{i\in I}$ with corresponding distinguished submersions $s_i\!:\!U_i\!\rightarrow\! \cT_i$. By compactness of $M$, we can assume our cover to be finite. The holonomy transformations (Haefliger's cocycles) $h_{ij}\!:\!s_i(U_i\cap U_j)\!\rightarrow\! s_j(U_i\cap U_j)$ induced by the transition maps on $M$ allow to define a manifold structure on the disjoint union $\cT\!=\!\bigsqcup_{i\in I} \cT_i$. This defines a complete transversal as required.\smallskip\\
\indent Recall that a foliation $\cF$ is a \textbf{Riemannian foliation} if the normal bundle $\nu\cF$ is equipped with a transverse metric. Equivalently, $\cF$ is a Riemannian foliation if on a complete transversal there exists a Riemannian metric such that all Haefliger's cocycles are isometries. For simplicity, we say that $\cF$ is Riemannian.
\begin{Def}[Holonomy pseudogroup, Salem\cite{Sal1}] The \textbf{holonomy pseudogroup} $\mathit{Hol}$ on $\cT$ is the pseudogroup generated by Haefliger's cocycles $h_{ij}$.
\end{Def}
\begin{Prop}\label{HolRiemannian}
The foliation $\cF$ is Riemannian if there exists a $\mathit{Hol}$-invariant Riemannian metric on $\cT$.
\end{Prop}
\subsection{Differential operators on foliated manifold}\ \\ \\
\indent The basic concepts of operator theory can be found in \cite{Kad,Kad1,BGV}. Let $P:C^\infty(M)\rightarrow C^\infty(M)$ be a differential operator of order $m\!\geq \! 2$ acting on the complex-valued smooth functions $C^\infty(M)$ of $M$. 
\begin{Def}\label{TEBO}
We say that $P$ is \textbf{basic} if at each point of $M$, there exists a foliated chart $U$ with coordinates $(x,y)\!=\!(x_1,\cdots,x_p,y_1,\cdots,y_q)$ such that $P$ can be written on $U$ as follows
\begin{equation}\label{eq:1.1}
\textstyle P\!=\!\sum_{|s|=0}^m a_s(y)\dfrac{\partial^{|s|}}{\partial y_1^{s_1}\cdots \partial y_q^{s_q}},
\end{equation}
where $s\!=\!(s_1,\cdots,s_q)\!\in\!\mathbb{N}^q$, $|s|\!=\!s_1\!+\!\cdots\!+\!s_q$ and $a_s$ are complex-valued smooth functions depending only on $y$.
\end{Def}
\begin{Rem}
The meaning of ``basic'' is different from \cite[Definition 2.7.2]{AAHR}. Indeed, El Kacimi's operator acts on global basic sections of a vector $\cF$-bundle, $\cF$ is supposed to be Riemannian.
\end{Rem}
\noindent In the sequel, let $P$ be a basic operator as in Definition \ref{TEBO}. In a foliated chart $U\!\simeq\! (x,y)$, a subspace defined by $dx_1\!=\!\cdots\!=\!dx_p\!=\!0$ is called \textbf{vertical subspace}. (The leaf directions defined by $dy_1\!=\!\cdots\!=\!dy_q\!=\!0$ are thought to be horizontal.) Vertical subspaces $\mathrm{T}_i$, $\mathrm{T}_j$ in local charts correspond to transversals $\cT_i$ and $\cT_j$, see Figure 1 below. 
\begin{Def}
A transversal $\cT$ is \textbf{compatible} with $P$ if at each point of $\cT$ there exists a foliated chart $U$ such that $U\!\cap\! \cT$ belongs to a vertical subspace.
\end{Def}
\begin{Rem}
In general, at a point of $M$, a compatible transversal is not unique. In an other foliated chart $\widetilde{U}\!\simeq\! (\tilde{x},\tilde{y})$, the expression of $P$ may also have the form $($\ref{eq:1.1}$)$ with an other expression. Notice that the vertical subspaces defined by $d\tilde{x}_1\!=\!\cdots\!=\!d\tilde{x}_p\!=\!0$ and $dx_1\!=\!\cdots\!=\!dx_p\!=\!0$ do not coincide. This gives that $\cT_i$ and $\widetilde{\cT}_i$ are both compatible transversals at $z$, see Figure 2 below. 
\end{Rem}
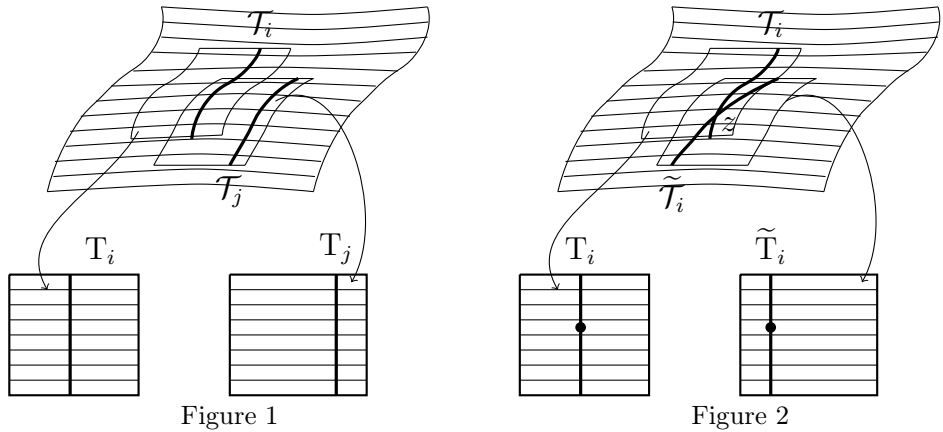
\begin{center}
\begin{tikzpicture}[domain=0:5]
\draw plot[smooth, tension=1] coordinates {(1.2,0.4) (1.5,1.15) (2,1.75) (2.6,2.3) (2.7,2.85)};
\draw plot[smooth] coordinates {(2.7,2.85) (4.7,2.75) (6.2,2.85)};
\draw plot[smooth, tension=1] coordinates {(4.7,0.4) (5,1.15) (5.5,1.75) (6.1,2.3) (6.2,2.85)};
\draw plot[smooth] coordinates {(1.2,0.4) (3.2,0.5) (4.7,0.4)};
\draw plot[smooth] coordinates {(1.28,0.6) (3.18,0.7) (4.72,0.6)};
\draw plot[smooth] coordinates {(1.3,0.8) (3.2,0.9) (4.8,0.8)};
\draw plot[smooth] coordinates {(1.44,1) (3.34,1.1) (4.92,1)};
\draw plot[smooth] coordinates {(1.5,1.2) (3.4,1.3) (5.05,1.2)};
\draw plot[smooth] coordinates {(1.64,1.4) (3.54,1.5) (5.18,1.4)};
\draw plot[smooth] coordinates {(1.84,1.6) (3.74,1.7) (5.32,1.6)};
\draw plot[smooth] coordinates {(2.04,1.8) (3.94,1.85) (5.52,1.8)};
\draw plot[smooth] coordinates {(2.33,2) (4.22,2.05) (5.81,2)};
\draw plot[smooth] coordinates {(2.58,2.25) (4.52,2.20) (6.05,2.25)};
\draw plot[smooth] coordinates {(2.7,2.45) (4.7,2.35) (6.2,2.45)};
\draw plot[smooth] coordinates {(2.71,2.65) (4.71,2.55) (6.21,2.65)};
\draw plot[smooth, tension=1] coordinates {(2.3,1.1) (2.5,1.6) (3,2) (3.2,2.3)};
\draw plot[smooth, tension=1] coordinates {(2.3,1.1) (3.5,1.15)};
\draw plot[smooth, tension=1] coordinates {(3.5,1.15) (3.7,1.6) (4.2,2) (4.4,2.3)};
\draw plot[smooth, tension=1] coordinates {(3.2,2.3) (4.4,2.3)};
\draw plot[smooth, tension=1] coordinates {(2.6,0.75) (2.8,1.1) (3.1,1.6) (3.5,1.9)};
\draw plot[smooth, tension=1] coordinates {(2.6,0.75) (3.8,0.75)};
\draw plot[smooth, tension=1] coordinates {(3.8,0.75) (4,1.1) (4.3,1.6) (4.7,1.9)};
\draw plot[smooth, tension=1] coordinates {(3.5,1.9) (4.7,1.9)};
\draw[thick] (0.7,-0.7) -- (0.7,-2.3) -- (2.4,-2.3) -- (2.4,-0.7) -- (0.7,-0.7);
\draw (0.7,-0.9) -- (2.4,-0.9);
\draw (0.7,-1.1) -- (2.4,-1.1);
\draw (0.7,-1.3) -- (2.4,-1.3);
\draw (0.7,-1.5) -- (2.4,-1.5);
\draw (0.7,-1.7) -- (2.4,-1.7);
\draw (0.7,-1.9) -- (2.4,-1.9);
\draw (0.7,-2.1) -- (2.4,-2.1);
\draw[thick] (3.6,-0.7) -- (3.6,-2.3) -- (5.4,-2.3) -- (5.4,-0.7) -- (3.6,-0.7);
\draw (3.6,-0.9) -- (5.4,-0.9);
\draw (3.6,-1.1) -- (5.4,-1.1);
\draw (3.6,-1.3) -- (5.4,-1.3);
\draw (3.6,-1.5) -- (5.4,-1.5);
\draw (3.6,-1.7) -- (5.4,-1.7);
\draw (3.6,-1.9) -- (5.4,-1.9);
\draw (3.6,-2.1) -- (5.4,-2.1);
\draw[very thick] (1.5,-0.7) -- (1.5,-2.3);
\draw[very thick] (5,-0.7) -- (5,-2.3);
\node[above] at (1.9,-0.7) {$\mathrm{T}_i$};
\node[above] at (5,-0.7) {$\mathrm{T}_j$};
\draw[->] (2.4,1.2) to [out=240,in=120] (1.2,-0.9);
\draw[->] (4.2,1.6) to [out=30, in=60] (5.2,-0.8);
\draw[very thick] plot[smooth, tension=1] coordinates {(3.1,1.1) (3.3,1.6) (3.8,2) (4,2.3)};
\draw[very thick] plot[smooth, tension=1] coordinates {(3.6,0.75) (3.8,1.1) (4.1,1.6) (4.5,1.9)};
\node[above] at (4,2.3) {$\mathcal{T}_i$};
\node[below] at (3.6,0.75) {$\mathcal{T}_j$};
\node at (3.6,-2.6) {\footnotesize Figure 1};
\end{tikzpicture}
\hspace{2em}
\begin{tikzpicture}[domain=0:4]
\draw plot[smooth, tension=1] coordinates {(1.2,0.4) (1.5,1.15) (2,1.75) (2.6,2.3) (2.7,2.85)};
\draw plot[smooth] coordinates {(2.7,2.85) (4.7,2.75) (6.2,2.85)};
\draw plot[smooth, tension=1] coordinates {(4.7,0.4) (5,1.15) (5.5,1.75) (6.1,2.3) (6.2,2.85)};
\draw plot[smooth] coordinates {(1.2,0.4) (3.2,0.5) (4.7,0.4)};
\draw plot[smooth] coordinates {(1.28,0.6) (3.18,0.7) (4.72,0.6)};
\draw plot[smooth] coordinates {(1.3,0.8) (3.2,0.9) (4.8,0.8)};
\draw plot[smooth] coordinates {(1.44,1) (3.34,1.1) (4.92,1)};
\draw plot[smooth] coordinates {(1.5,1.2) (3.4,1.3) (5.05,1.2)};
\draw plot[smooth] coordinates {(1.64,1.4) (3.54,1.5) (5.18,1.4)};
\draw plot[smooth] coordinates {(1.84,1.6) (3.74,1.7) (5.32,1.6)};
\draw plot[smooth] coordinates {(2.04,1.8) (3.94,1.85) (5.52,1.8)};
\draw plot[smooth] coordinates {(2.33,2) (4.22,2.05) (5.81,2)};
\draw plot[smooth] coordinates {(2.58,2.25) (4.52,2.20) (6.05,2.25)};
\draw plot[smooth] coordinates {(2.7,2.45) (4.7,2.35) (6.2,2.45)};
\draw plot[smooth] coordinates {(2.71,2.65) (4.71,2.55) (6.21,2.65)};
\draw plot[smooth, tension=1] coordinates {(2.3,1.1) (2.5,1.6) (3,2) (3.2,2.3)};
\draw plot[smooth, tension=1] coordinates {(2.3,1.1) (3.5,1.15)};
\draw plot[smooth, tension=1] coordinates {(3.5,1.15) (3.7,1.6) (4.2,2) (4.4,2.3)};
\draw plot[smooth, tension=1] coordinates {(3.2,2.3) (4.4,2.3)};
\draw plot[smooth, tension=1] coordinates {(2.5,0.75) (2.7,1.1) (3,1.6) (3.4,1.9)};
\draw plot[smooth, tension=1] coordinates {(2.5,0.75) (3.7,0.75)};
\draw plot[smooth, tension=1] coordinates {(3.7,0.75) (3.9,1.1) (4.2,1.6) (4.6,1.9)};
\draw plot[smooth, tension=1] coordinates {(3.4,1.9) (4.6,1.9)};
\draw[thick] (0.7,-0.7) -- (0.7,-2.3) -- (2.4,-2.3) -- (2.4,-0.7) -- (0.7,-0.7);
\draw (0.7,-0.9) -- (2.4,-0.9);
\draw (0.7,-1.1) -- (2.4,-1.1);
\draw (0.7,-1.3) -- (2.4,-1.3);
\draw (0.7,-1.5) -- (2.4,-1.5);
\draw (0.7,-1.7) -- (2.4,-1.7);
\draw (0.7,-1.9) -- (2.4,-1.9);
\draw (0.7,-2.1) -- (2.4,-2.1);
\draw[thick] (3.6,-0.7) -- (3.6,-2.3) -- (5.4,-2.3) -- (5.4,-0.7) -- (3.6,-0.7);
\draw (3.6,-0.9) -- (5.4,-0.9);
\draw (3.6,-1.1) -- (5.4,-1.1);
\draw (3.6,-1.3) -- (5.4,-1.3);
\draw (3.6,-1.5) -- (5.4,-1.5);
\draw (3.6,-1.7) -- (5.4,-1.7);
\draw (3.6,-1.9) -- (5.4,-1.9);
\draw (3.6,-2.1) -- (5.4,-2.1);
\draw[->] (2.4,1.2) to [out=240,in=120] (1.2,-0.9);
\draw[->] (4.2,1.6) to [out=30, in=60] (5.2,-0.8);

\draw[very thick] (1.5,-0.7) -- (1.5,-2.3);
\draw[very thick] (4,-0.7) -- (4,-2.3);
\node[above] at (1.5,-0.7) {$\mathrm{T}_i$};
\node[above] at (4,-0.7) {$\widetilde{\mathrm{T}}_i$};
\draw[very thick] plot[smooth, tension=1] coordinates {(3.2,1.1) (3.4,1.6) (3.9,2) (4.1,2.3)};
\draw[very thick] plot[smooth, tension=1] coordinates {(2.7,0.75) (2.9,1) (3.4,1.5) (4.1,1.9)};
\node[above] at (4,2.3) {$\mathcal{T}_i$};
\node[below] at (2.7,0.75) {$\widetilde{\mathcal{T}}_i$};
\fill (1.5,-1.4) circle (0.07);
\fill (4,-1.4) circle (0.07);
\node[right] at (3.2,1.3) {$z$};
\node at (3.6,-2.6) {\footnotesize Figure 2};
\end{tikzpicture}
\end{center}

\begin{Prop}\label{ConstructionTransversal}
There exists a complete transversal $\cT\!=\!\bigsqcup_{i\in I} \cT_i$ of $\cF$ satisfying
\begin{center}
\noindent\fcolorbox{black}{white}{
\begin{minipage}{0.7\textwidth}
\noindent\textbf{Assumption:}
\begin{itemize}
\item[$(1)$] $\cT$ has a finite number of open connected components $\cT_i$;
\item[$(2)$] Each $\{\cT_i\}$ is a relatively compact subset of $\{\cT_i^\prime\}$;
\item[$(3)$] Each $\cT_i^\prime$ being compatible to $P$ locates in a foliated chart where $P$ has the form $(1)$.
\end{itemize}
\end{minipage}
}
\end{center}
\end{Prop}
\Preuve For every $z\!\in\! M$, there exists a foliated chart $U_z^\prime\!\simeq\!\mathbb{R}^q\!\times\!\mathbb{R}^q$ such that $P$ has the form (\ref{eq:1.1}) and $z\!\simeq\!0$. Let $U_z\!\simeq\!\mathbb{R}^p\!\times\! B(0,1)$. Then, the cover $\{U_z\}_{z\in M}$ admits a finite cover $\{U_i\}_{i\in I}$. For each $i\!\in\! I$, we define $\cT_i\!\simeq\!\{0\}\!\times\!B(0,1)$ and $\cT_i^\prime\!\simeq\!\{0\}\!\times\!\mathbb{R}^q$. Then, $\cT\!=\!\bigsqcup_{i\in I} \cT_i$ is a complete transversal satisfying \textbf{Assumption}.\eb
In the following, $\cT$ always satisfies \textbf{Assumption}. Clearly, $\cT$ is compatible to $P$ and $\cT$ is a relatively compact subset of $\cT^\prime\!=\!\bigsqcup_{i\in I}\cT_i^\prime$.\smallskip\\
\indent In \cite{Hae}, Haefliger proved that when $M$ is compact, the holonomy pseudogroup \textit{Hol} of $\cF$ satisfies the property of compact generation. For a compactly generated pseudogroup of local isometries of a Riemannian manifold, Salem proposed a version of Molino's theory in \cite{Sal}, \cite{Sal1}, also \cite{AM}.
\begin{Def}\label{CompactlyGenerated}
A pseudogroup $\bH$ on $\bT$ is said to be \textbf{compactly generated} if
\begin{itemize}
\item[$\bullet$] $\exists$ a relatively compact open subset ${\mathrm{\bf U}}\!\subset\!\bT$ meeting all $\bH$-orbits;
\item[$\bullet$] $\exists$ a finite set $S\!=\!\{h_1,\cdots,h_n\}\!\subset\!\bH\vert_{\mathrm{\bf U}}$ that generates $\bH\vert_{\mathrm{\bf U}}$;
\item[$\bullet$] each $h_i$ is the restriction of some $\widetilde{h_i}\!\in\!\bH$ with $\overline{\operatorname{dom}h_i}\!\subset\! \operatorname{dom}\widetilde{h_i}$.
\end{itemize}
\end{Def}
\begin{Prop}\label{HolFinite}
The holonomy pseudogroup \textit{Hol} on $\cT^\prime$ is compactly generated. 
\end{Prop}
\Preuve Let $\bT\!=\!\cT^\prime,\,\mathrm{\bf U}\!\!=\!\!\cT$ in Definition \ref{CompactlyGenerated}, $\cT$ meets all $\Hol$-orbits because $\Hol$-orbits are exactly the leaves. Since $I$ is a finite set, the Haefliger's cocycle set $\{h_{ij}\mid (i,j)\in I^2\}$ is also finite. By \textbf{Assumption} (2), $h_{ij}$ is a restriction of Haefliger's cocycle $h_{ij}^\prime$ defined on $\cT^\prime$. \eb
Restricted to $\cT$, $P:C^\infty(\cT)\rightarrow C^\infty(\cT)$ is a differential operator. Till the end of this section, $P$ refers to the restricted operator. 
\begin{Def}
A local transformation $\psi$ on $\cT$ commutes with $P$ if for every $f\!\in\! C^\infty(\operatorname{im}\psi)$ the relation
\begin{equation}\label{Commutator}
\psi^* (Pf)\!=\!P(\psi^*f)
\end{equation}
holds on $\operatorname{dom}\psi$ where $\psi^*f\!\!=\!\!f\!\circ\!\psi$.
\end{Def}
\begin{Prop}\label{GPHol}
The collection $\bG(P)$ of all local transformations on $\cT$ commuting with $P$ is a pseudogroup, called the \textbf{commuting pseudogroup} of $P$. The holonomy pseudogroup on $\cT$ is a sub-pseudogroup of $\bG(P)$.
\end{Prop}
\Preuve It is easy to check that $\bG(P)$ is a pseudogroup. The local transformations of holonomy commute with $P$ because in expression $($\ref{eq:1.1}$)$, the operator $P$ only depends on $y$.\eb
\begin{Prop}\label{Nozeropart}
If a local transformation $\psi$ belongs to $\bG(P)$, then $\psi$ commutes with the part $P_{\geq 1}$ of order $\geq\! 1$.
\end{Prop}
\Preuve For each $\psi\!\in\!\bG(P)$, its domain may have many connected components. By Definition \ref{pseudogroup} (5), $\psi$ can be regarded as the combination of its restrictions on each connected components. Thus, we may assume that $\operatorname{dom}\psi$ and $\operatorname{im}\psi$ locate respectively in some $\cT_i$ and $\cT_j$. By \textbf{Assumption} (3), $P$ is written on $\operatorname{dom}\psi$ and $\operatorname{im}\psi$ by
$$
\textstyle P\!=\!a_0(y)\!+\!\sum_{|s|=1}^m a_s(y)\dfrac{\partial^{|s|}}{\partial y_1^{s_1}\cdots \partial y_q^{s_q}};\ \ 
P\!=\!b_0(y^\prime)\!+\!\sum_{|s|=1}^m b_s(y^\prime)\dfrac{\partial^{|s|}}{\partial {y^\prime_1}^{s_1}\cdots \partial {y^\prime_q}^{s_q}};
$$
where $y,y^\prime$ are coordinates on $\operatorname{dom}\psi,\operatorname{im}\psi$, respectively. And $y^\prime=\psi(y)$.\\
In relation (\ref{Commutator}), take the constant function $f\!=\!1$ on $\operatorname{im}\psi$, then we get
$b_0\big(\psi(y)\big)\!=\!a_0(y)$. Thus the zero-order parts on two sides cancel each other. \eb
\section{Commuting Pseudogroup}
Let the foliated manifold $(M,\cF)$ be compact. In this section, we aim to obtain the properties of equicontinuity and quasi-analyticity of the commuting pseudogroup and summarize to the following theorem:
\begin{Thm}\label{Thm1}
Let $(M,\cF)$ be a compact foliated manifold and $P\!:\!C^\infty(M)\!\rightarrow\! C^\infty(M)$ be a basic differential operator, $\cT$ be a complete transversal satisfying \textrm{\bf Assumption} and $\bG(P)$ be the commuting pseudogroup.
\begin{enumerate}
\item If $P$ is transversely elliptic, $\bG(P)$ is equicontinous with respect to the induced Riemannian distance;
\item If $P$ is transversely elliptic and analytic with triangular $1$-part at $z\!\in\!\cT$, $\bG(P)$ is quasi-analytic around $z$;
\item A sub-pseudogroup $\cH$ of $\bG(P)$ is quasi-analytic if its restriction to an open subset of $\cT$ meeting all $\cH$-orbits is quasi-analytic.
\end{enumerate}
\end{Thm}
\subsection{Equicontinuity of $\bG(P)$}\ \\ \\
\indent Let us recall the notion of transverse ellipticity of a differential operator. Denote by $(\nu\cF)^*\!=\!\big\{(z,\xi)\in T^*M\mid \forall\, v\in T\cF\vert_z, \xi(v)=0\big\}$.
\begin{Def}
We say that $P$ is \textbf{transversely elliptic} if for every $z\!\in\! M$, in expression $($\ref{eq:1.1}$)$, the principal symbol
$$
\textstyle\sigma(P)_z(\xi)\!=\!\sum_{|s|=m}a_s(y)\xi_1^{s_1}\cdots \xi_q^{s_q}
$$
is nonzero for every non-zero transverse covector $\xi\!\in\!(\nu\cF)^*\vert_z$.
\end{Def}

\indent In this section, let $P$ be a basic transversely elliptic operator of order $m\!\geq\!2$. The restricted operator $P\!:\!C^\infty(\cT)\!\rightarrow\! C^\infty(\cT)$ is  elliptic. Let $\xi$ be a cotangent vector in $T^*_z\cT$. Take a function $f$ on $\operatorname{im}\psi$ such that $f\vert_z\!\!=\!\!0$ and $(df)\vert_{z}\!\!=\!\!\xi$ where $\vert$ representing evaluation. The principal symbol of $P$ at $(z,\xi)$ is also given by
$$
\textstyle\sigma(P)_z(\xi)\!=\!\dfrac{1}{m!}(Pf^m)\vert_z.
$$
\begin{Lemme}\label{Symbol}
For every $\psi\!\in\!\bG(P)$, $z\!\in\!\operatorname{im}\psi$, $\xi\!\in\! T^*_z\cT$, we have
$$
\textstyle\sigma(P)_{\psi^{-1}(z)}(\psi^*\xi)\!=\!\sigma(P)_z(\xi).
$$
\end{Lemme}
\Preuve  It is easy to see $(\psi^*f)\vert_{\psi^{-1}(z)}\!=\!0$ and $(d\psi^*f)\vert_{\psi^{-1}(z)}\!=\!\psi^*\xi$. A brief calculation:
$$
\textstyle \sigma(P)_{\psi^{-1}(z)}(\psi^*\xi)\!=\!\dfrac{1}{m!}\big(P(\psi^*f)^m\big)\vert_{\psi^{-1}(z)}\!=\!\dfrac{1}{m!}\big(\psi^*Pf^m\big)\vert_{\psi^{-1}(z)}
\!=\!\dfrac{1}{m!}(Pf^m)\vert_z\!=\!\sigma(P)_z(\xi).
$$ 
\\[-2em]\eb
\medskip
Fix a Riemannian metric on $M$, which induces a Riemannian metric on $\cT$. Let $|\!|\cdot|\!|$ denote the induced norm on the tangent bundle $T\cT$. For every $\psi\!\in\!\bG(P)$, $z\!\in\!\operatorname{dom}\psi$, the tangent map $T_z\psi$ is normed by $|\!|T_z\psi|\!|\!=\!\mathrm{sup}_{|\!|v|\!|=1} |\!|T_z\psi(v)|\!|$. 
\begin{Prop}\label{Bounds}
There exist uniform bounds $\lambda,\mu>0$ such that $\forall\, \psi\!\in\!\bG(P)$, $\forall\, z\!\in\!\operatorname{dom}\psi$, 
$$
\lambda \!\leq\! |\!|T_z\psi|\!|\!\leq\! \mu.
$$
\end{Prop}
\Preuve  \underline{\textbf{Upper-boundness}.}\\
Suppose that there exist a sequence $\psi_n\!\in\!\bG(P)$, a sequence $z_n\!\in\!\cT$ and $v_n\!\in\! T^1_{z_n}\cT$ with $|\!|v_n|\!|\!=\!1$ such that $|\!|T_{z_n}{\psi_n}(v_n)|\!|\!\rightarrow\! +\infty$. Under \textbf{Assumption}, $\overline{\cT}$ is compact. Then we assume $z_n\!\!\rightarrow\!\! z\!\in\!\overline{\cT_i}$ and $w_n\!=\!\psi_n(z_n)\!\rightarrow\! w\!\in\!\overline{\cT_j}$ for some $i,j\!\in\! I$. \\
For $n$ large enough, $z_n$ (resp. $w_n$) belongs to $\cT_i^\prime$ (resp. $\cT_j^\prime$) with coordinates $u$ (resp. $v$). By denoting simply ${\psi_n}_*\!=\!T_{z_n}{\psi_n}$, one has
$$
\textstyle{\psi_n}_* \begin{bmatrix}\begin{smallmatrix}\frac{\partial}{\partial u_1}\vert_{z_n}& \cdots & \frac{\partial}{\partial u_q}\vert_{z_n}\end{smallmatrix}\end{bmatrix}^T\!=\!A_n \begin{bmatrix}\begin{smallmatrix}\frac{\partial}{\partial v_1}\vert_{w_n}&\cdots & \frac{\partial}{\partial v_q}\vert_{w_n}\end{smallmatrix}\end{bmatrix}^T,
$$
where $[\ ]^T$ is the transposition of matrix, the real matrix $A_n\!=\!(a_{kl}^n)$ is the Jacobian $J(\psi_n)$ of $\psi_n$ evaluated at $z_n$. By duality, we have
$$
\psi^*_n\begin{bmatrix}\begin{smallmatrix}dv_1\vert_{w_n}&\cdots&dv_q\vert_{w_n}\end{smallmatrix}\end{bmatrix}\!=\!\begin{bmatrix}\begin{smallmatrix}du_1\vert_{z_n}&\cdots&du_q\vert_{z_n} \end{smallmatrix}\end{bmatrix}A_n.
$$
Denote simply by $\sigma(P)(\xi)$ for $\sigma(P)_z(\xi)$ in this proof.
On one hand, by Lemma \ref{Symbol}, 
$$
\textstyle\sigma(P)\big(\sum_{i=1}^q a_{il}^n du_i\vert_{z_n}\big)\!=\!\sigma(P)\big(\psi^*_n dv_l\vert_{w_n}\big)\!=\!\sigma(P)(dv_l\vert_{w_n})\stackrel{n\rightarrow \infty}{\longrightarrow}\sigma(P)(dv_l\vert_w)\!\neq\! 0.
$$ 
If the $kl$-component $a_{kl}^n$ of $A_n$ is unbounded and assume $|a_{kl}^n|\!\rightarrow\! +\infty$, on the other hand 
$$
\textstyle\sigma(P)\big(\sum_{i=1}^q a_{il}^n du_i\vert_{z_n}\big)
\!=\!(a_{kl}^n)^{m} \sigma(P)\Big(du_k\vert_{z_n}\!+\!
\sum_{i\neq k} \frac{a_{il}^n}{a_{kl}^n}du_i\vert_{z_n}\Big).
$$
Thus we have
$$
\textstyle\sigma(P)\big(du_k\vert_{z_n}\!+\!\sum_{i\neq k} \frac{a_{il}^n}{a_{kl}^n}du_i\vert_{z_n}\big)\!\rightarrow 0.
$$ 
By the smoothness and ellipticity of $P$, $\textstyle du_k\vert_{z_n}\!\!+\!\!
\sum_{i\neq k} \frac{a_{il}^n}{a_{kl}^n}du_i\vert_{z_n}\!\!\rightarrow\! 0,$ which contradicts $du_k\vert_{z_n}\!\!\rightarrow \!du_k\vert_z\!\neq\! 0$. Therefore, each component of $A_n$ is bounded.\\
One has a subsequence, denoted also by $\{A_n\}$ convergent to a real matrix $A$. On the other hand, the compactness of the unit bundle $(T^1\cT_i^\prime)\vert_{\overline{\cT_i}}$ implies a subsequence $\{v_n\}$ convergent to $v_\infty\!=\!\sum_{i=1}^q\alpha_i\frac{\partial}{\partial v_i}\vert_z$. Then ${\psi_n}_*(v_n)$ converges to
$$
\begin{bmatrix}\begin{smallmatrix}\alpha_1&\cdots&\alpha_q\end{smallmatrix}\end{bmatrix}A\begin{bmatrix}\begin{smallmatrix}
\frac{\partial}{\partial v_1}\vert_{w}&  \cdots& \frac{\partial}{\partial v_q}\vert_{w} \end{smallmatrix}\end{bmatrix}^{T},
$$
which contradicts the divergence of $|\!|{\varphi_n}_*(v_n)|\!|$.\\
\underline{\textbf{Lower Boundness}}\\ 
\noindent Similarly, if $|\!|{\varphi_n}_*(v_n)|\!|\!\rightarrow\! 0$, all subsequences converge to zero. Hence, 
$$
\begin{bmatrix}\alpha_1&\cdots&\alpha_q\end{bmatrix}A\!=\!0.
$$ 
Since $(\alpha_1,\cdots,\alpha_q)\!\neq\!0$ ( $|\!|v_\infty|\!|\!\!=\!\!1$), $A$ has zero determinant. There is $(\beta_1,\cdots, \beta_q)\!\neq\! 0$ such that $A\begin{bmatrix}\beta_1& \cdots & \beta_q \end{bmatrix}^T\!=\!0$. The sequence of covectors $\xi_n\!=\!\sum_{i=1}^q \beta_i dv_i\vert_{w_n}$ converges to the non-zero covector $\xi\!=\!\sum_{i=1}^q \beta_i dv_i\vert_{w}$. Thus, $\sigma(P)(\xi_n)\rightarrow\sigma(P)(\xi)\!\neq\! 0$.\smallskip\\
But in the relation $\sigma(P)(\varphi_n^*\xi_n)\!=\!\sigma(P)(\xi_n)$, the left side converges to $0$:
$$
\textstyle\psi_n^*(\xi_n)\!=\!\begin{bmatrix}du_1\vert_{z_n}&\cdots&du_q\vert_{z_n} \end{bmatrix}A_n\begin{bmatrix}\beta_1\\ \vdots \\ \beta_q\end{bmatrix}\rightarrow \begin{bmatrix} du_1\vert_{z}&\cdots&du_q\vert_{z}\end{bmatrix}A\begin{bmatrix}\beta_1\\ \vdots \\ \beta_q\end{bmatrix}\!=\!0
$$
\vspace{-1.5em} \eb
\medskip
Let us introduce the notion of equicontinuity of pseudogroup. Here, we use a notion slightly stronger than the one used in \cite[Definition 2.29]{AG}. 
\begin{Def}[Equicontinuity] Let $\bT\!=\!\bigsqcup_{i\in I} \bT_i$ be a smooth manifold. A pseudogroup $\bH$ on $\bT$ is said to be \textit{equicontinuous} if there is some metric $d_i$ on each component $\bT_i$ for all $i\!\in\! I$ such that $\forall\,\epsilon\!>\!0$, $\exists\,\delta(\epsilon)$ so that 
$$
d_i(z,w)\!<\!\delta(\epsilon)\Longrightarrow d_j\big(h(z),h(w)\big)\!<\!\epsilon
$$ 
for all $h\!\in\!\bH$, $i,j\!\in\! I$ and $z,w\!\in\!\bT_i\cap h^{-1}(\bT_j\cap \operatorname{im}h)$.
\end{Def}
In the present work, each component $\cT_i$ of $\cT$ has an induced Riemannian metric. Consider the Riemannian distance defined by $d_i(z,w)\!=\!\inf_{c} \int_0^1 |\!|c^\prime(t)|\!|dt$ where $c$ runs over all paths from $z$ to $w$.
\begin{Prop}\label{Equicontinuity}
The pseudogroup $\bG(P)$ is equicontinous with respect to the Riemannian distance.
\end{Prop}
\Preuve By Proposition \ref{Bounds}, we have a uniform upper bound $\mu$. For every $\epsilon>0$, take $\delta(\epsilon)\!=\!\frac{1}{\mu}\epsilon$. For all $\psi\!\in\! \bG(P)$, $i,j\!\in\! I$ and $z,w\!\in\!\cT_i\cap \psi^{-1}(\cT_j\cap\operatorname{im}\psi)$ with $d(z,w)\!<\!\delta(\epsilon)$, for each path $c$ from $z$ to $w$ located in $\operatorname{dom}\psi$, $\varphi\!\circ\! c$ is a path from $\psi(z)$ to $\psi(w)$ located in $\operatorname{im}\psi$. Thus,
$$
\textstyle d\big(\psi(z),\psi(w)\big)\!\leq\! \int_0^1 |\!|(\psi\!\circ\! c)^\prime(t)|\!|dt\!\leq\!\int_0^1 |\!|T_{c(t)}\psi|\!| |\!|c^\prime(t) |\!|dt\!\leq\! \mu \int_0^1 |\!|c^\prime(t)|\!|dt.
$$
Take the infimun for $c$ running over all paths from $z$ to $w$, therefore
$$
\textstyle d\big(\psi(z),\psi(w)\big)\!\leq\!\mu\inf_{c} \int_0^1 |\!|c^\prime(t)|\!|dt\!=\!\mu d(z,w)\!<\!\mu\delta(\epsilon)\!=\!\epsilon.
$$
\vspace{-1.5em}\eb 
We now present the following counter-example already studied in \cite{Liu1}, see also \cite[Section 4.1.1]{Liu}.
\begin{CExa} 
Let $M\!=\!S^1\times S^1\times S^1$ with coordinates denoted by $(x,y,z)\!\in\! M$. Let $f\!\in\! C^\infty(S^1)$ be a non-constant real function. The foliation $\cF$ is generated by $V\!=\!\frac{\partial}{\partial x}\!+\!f(y)\frac{\partial}{\partial z}$. Without losing generality, assume that at a point $m\!\!=\!\!(x_0,y_0,z_0)$, one has $f(y_0)\!=\!0$, $f^\prime(y_0)\!=\! 1$. Consider $\cT\!=\!\{x_0\}\times S^1\times S^1$ to be a $($complete$)$ transversal. The path $c_n:[0,1]\rightarrow M,\, t\mapsto (x_0e^{2n\pi it},y_0,z_0)$ located in the leaf of $\cF$ passing by $m$ determines a local diffeomorphism $\psi_n$. A brief calculation shows 
$$
\textstyle T_m{\psi_n}\begin{bmatrix}\begin{smallmatrix}\frac{\partial}{\partial y}\vert_{m}\\ \\ \frac{\partial}{\partial z}\vert_{m}\end{smallmatrix}\end{bmatrix}\!=\!\begin{bmatrix}\begin{smallmatrix}1&n\\ \\ 0&1\end{smallmatrix}\end{bmatrix}\begin{bmatrix}\begin{smallmatrix}\frac{\partial}{\partial y}\vert_{m}\\ \\ \frac{\partial}{\partial z}\vert_{m} \end{smallmatrix}\end{bmatrix}.
$$
Then $A_n\!=\!\begin{bmatrix}\begin{smallmatrix}1&n\\ \\ 0&1\end{smallmatrix}\end{bmatrix}$ and therefore the symmetric matrix $A_n^TA_n$ has a divergent eigenvalue $\lambda\!=\! \frac{1}{2}\big(n^2\!+\!2\!+\!n\sqrt{n^2\!+\!4}\,\big)$, so $|\!|T_m\psi_n|\!|\!\rightarrow\! \infty$. It follows that no basic transversely elliptic operator $P$ exists on $(M,\cF)$.
\end{CExa}

\subsection{Quasi-Analyticity of $\bG(P)$}\ \\ \\
\indent In this section, let $P\!:\!C^\infty(M)\!\rightarrow\! C^\infty(M)$ be a basic transversely elliptic operator of order $m\!\geq\!2$, $\cT$ be a complete transversal satisfying \textbf{Assumption}. \\
\indent Each local transformation $\psi$ of $\bG(P)$ may be thought to have a connected open domain. Let us introduce the notion of quasi-analyticity of pseudogroup. 
\begin{Def} \ \\[-1em]
\begin{itemize}
\item[$\bullet$] $($\cite{Hae1},\cite{AG}Definition 2.17$)$ A pseudogroup $\bH$ on $\bT$ is quasi-analytic if every $h\in\bH$ is the identity around some $z\in \operatorname{dom}h$ whenever $h$ is the identity on some open set whose closure contains $z$;
\item[$\bullet$] A pseudogroup $\bH$ is quasi-analytic around $z\in\bT$ if there exist an connected open neighborhood $\mathrm{\bf U}\subset \bT$ such that $\bH\vert_{\mathrm{\bf U}}$ is quasi-analytic.
\end{itemize}
\end{Def}
\begin{Prop}\label{SQAQA}
$\bH$ is quasi-analytic if and only if $\bH$ is quasi-analytic around every point $z\!\in\!\bT$.
\end{Prop}
\indent Suppose that $\psi\!\in\!\bG(P)$ is the identity on some open subset $V\!\subset\!\operatorname{dom}\psi$. Let $z\!\in\!\partial V$ and by smoothness of $\psi$, $\psi(z)\!=\!z$. Around $z$, $\psi$ can be regarded as a diffeomorphism of $\mathbb{R}^{q}$ which is the identity in the upper semi-plane. Then the $1$-jet $j^1_z\psi$ at $z$ is the identity. But $\psi$ is not necessarily the identity around $z$ because $\psi$ is not assumed to be isometric. (A local isometry is completely determined by its $1$-jet.) \\
\indent In fact, even if the $\infty$-jet $j^\infty_z\psi$ is the identity, $\psi$ is not necessarily the identity unless $\psi$ is proved to be analytic. Hence, a straight way to acquire the quasi-analyticity of $\bG(P)$ is to follow the analyticity of solution of a PDE system which $\psi$ satisfies. Under conditions in \cite{Mor}, $P$ is supposed to have analytic function coefficients. 
\begin{Def}
We say that $P$ is \textbf{analytic} at $z\!\in\! M$ if there exists a foliated chart at $z$ such that all function coefficients $a_s$ in expression $($\ref{eq:1.1}$)$ are analytic. $P$ is analytic on a subset of $M$ if $P$ is at each point. 
\end{Def}
By Proposition \ref{Nozeropart}, the zero-order part of $P$ does not induce any property of $\bG(P)$. Focusing on the order $1$-part, a supplementary condition is introduced.
\begin{Def}\label{Triangular1}
We say that $P$ is \textbf{with triangular $1$-part} at $z\!\in\! M$ if there exists a foliated chart $U$ of $z$ such that the following matrix   
$$
\dfrac{\partial(a_1,\cdots,a_q)}{\partial(y_1,\cdots,y_q)}=\begin{bmatrix}a_{11}&\cdots&a_{1q}\\ \vdots&\ddots&\vdots\\ a_{q1}&\cdots&a_{qq}
\end{bmatrix}
$$  
is triangular where the functions $(a_1,\cdots, a_q)$ are the coefficients of the order $1$-part of $P$
\begin{equation}\label{firstorderpart}
\textstyle P_1\!=\!a_1(y)\frac{\partial}{\partial y_1}\!+\!\cdots+a_q(y)\frac{\partial}{\partial y_q}
\end{equation} 
written in the coordinate $(x,y)$ of $U$. $P$ is said to be analytic with triangular $1$-part on a subset of $M$ if $P$ is at each point.
\end{Def}
Suppose that $P$ is analytic at a point $z$ belonging to the transversal $\cT$, then in a foliated chart $U\!\simeq\!(x,y)$, the vertical subspace $\mathrm{O}$ passing through $z$ is a compatible transversal. Notice that $\mathrm{O}$ is not necessarily a subset of $\cT$, whereas $\mathrm{O}$ and $\cT$ are both transversals at $z$, see Figure 2.\\
\indent Till the end of this section, $P$ is restricted to $O$. Without ambiguity, a point $w\!\in\! \mathrm{O}$ and its coordinates $y$ will not be distinguished. 
\begin{Lemme}
Each local transformation $\phi$ on $\mathrm{O}$ with $\operatorname{dom}\phi\cup\operatorname{im}\phi\!\subset\! \mathrm{O}$ commuting with $P$ is analytic if the determinant of the following matrix evaluated at every $w\!\in\! \operatorname{dom}\phi$, every non-zero covector $\eta\!\in\! T^*_w\mathrm{O}$ is non-zero: 
\begin{equation}\label{CharPoly}
\begin{bmatrix}
\sigma(P)&-\phi^*a_{12}&\cdots&-\phi^*a_{1q}\\
-\phi^*a_{21}&\sigma(P)&\cdots&-\phi^*a_{2q}\\ 
\vdots&\vdots&\ddots&\vdots\\
-\phi^*a_{q1}&-\phi^*a_{q2}&\cdots&\sigma(P)
\end{bmatrix},
\end{equation}   
where $a_{ij}$ are as in Definition \ref{Triangular1}.
\end{Lemme}
\Preuve Under coordinate $y$, let $\phi$ be written by $(\phi^1,\cdots,\phi^q)$. In relation (\ref{Commutator}), take $f(y)\!=\!y_k$, we have
$$
\phi^*a_k\!=\!\phi^*(Pf)\!=\!P(\phi^*f)\!=\!P\phi^k.
$$
Thus, the functions $\phi^1,\cdots,\phi^q$ satisfy the following non-linear PDE system on $\operatorname{dom}\phi$:
\begin{equation}\label{SystemPDE}
\begin{array}{l}
L_k(y):\ P\phi^k\!-\!a_k(\phi^1,\cdots,\phi^q)\!=\!0, \hspace{0.5cm} k\!=\!1,\cdots,q.
\end{array}
\end{equation}
The first variation of (\ref{SystemPDE}) (see \cite[Equation (1.3)]{Mor} for definition) along the function vector $v\!=\!(v^1,\cdots,v^q)$ at $\phi$ for $l\!=\!1,\cdots,q$ is
$$
\left\{
\begin{array}{l}
L_{kl}v_l\!=\frac{d}{d\lambda}\vert_{\lambda=0}\big[ P\phi^k\!-\!a_k(\phi^1,\cdots,\phi^l\!+\!\lambda v^l,\cdots,\phi^q)\big]\!=\!-\big(a_{kl}\!\circ\!\phi\big) v^l\ \  \ \text{if}\ l\!\neq\! k;\\ \\[-1em]
L_{kk}v_k\!=\!\frac{d}{d\lambda}\vert_{\lambda=0}\big[ P(\phi^k\!+\!\lambda v^k)\!-\!a_k(\phi^1,\cdots,\phi^k\!+\!\lambda v^k,\cdots,\phi^q)\big]\!=\!(Pv^k)\!-\!\big(a_{kk}\!\circ\!\phi\big) v^k.
\end{array}
\right.
$$
Let $L^0_{kl}$ be the top-order term of $L_{kl}$. The characteristic polynomial $L^0_{kl}(w,\eta)$ of $L^0_{kl}$ with respect to solution $\phi$ is 
$$
\left\{
\begin{array}{l}
L^0_{kl}(w,\eta)=-(\phi^* a_{kl})(w)\ \  \ \text{if}\ l\!\neq\! k;\\ \\[-1em]
L^0_{kk}(w,\eta)=\sigma(P)_w(\eta)
\end{array}
\right.
$$ 
Followed by \cite[Equation (1.4)]{Mor} and the analyticity of the solutions in the interior obtained by \cite[Theorem C]{Mor}, $\phi$ is analytic in the interior of its domain if the determinant
$$
\vert L^0_{kl}(w,\eta)\vert\neq 0
$$
for any $w\in \operatorname{dom}\phi$ and any non-zero $\eta\in T_w^*\mathrm{O}$. Since $\operatorname{dom}\phi$ is open and connected, $\phi$ is analytic. (The system $($\ref{SystemPDE}$)$ is said to be elliptic with respect to solution $\phi$.)\eb
The condition of triangular $1$-part admits the analyticity of all local transformations on $\mathrm{O}$ commuting with $P$, then induces the quasi-analyticity of $\bG(P)$ around a point $z\in\cT$.
\begin{Prop}\label{AtoQA}
The commuting pseudogroup $\bG(P)$ is quasi-analytic around $z\!\in\!\cT$ if $P$ is analytic at $z$ with triangular $1$-part.
\end{Prop}
\Preuve Under the condition of triangular $1$-part, the matrix $($\ref{CharPoly}$)$ reduces to a triangular matrix. Its determinant $[\sigma(P)_w(\eta)]^q$ is non-zero for every $(w,\eta)$ by the ellipticity of $P$ on $\mathrm{O}$. Hence, every $\phi$ with $\operatorname{dom}\phi\cup\operatorname{im}\phi\subset \mathrm{O}$ commuting with $P$ is analytic. \\
The holonomy in local chart $U$ defines a diffeomorphism $h\!:\!\mathrm{O}\!\rightarrow\! h(\mathrm{O})\!\subset\!\cT$. For $\psi\!\in\!\bG(P)\vert_{h(\mathrm{O})}$, the local transformation $\tilde{\psi}\!=\!h^{-1}\!\circ\! \psi\! \circ\! h$ on $\mathrm{O}$ commutes with $P$ because $P$ commutes with holonomy by Proposition \ref{GPHol}. Thus, if $\psi$ is identity on some open subset of $h(\mathrm{O})$, $\tilde{\psi}$ is identity on some open subset of $\mathrm{O}$. The analyticity of $\tilde{\psi}$ gives 
$$
\tilde{\psi}\!=\!\mathrm{id}_{\mathrm{O}}.$$ 
Therefore, $\psi\!=\!\mathrm{id}_{h(\mathrm{O})}$. \eb
In order to simplify notations, we shall denote by $\psi\psi^\prime$ the composition $\psi\!\circ\!\psi^\prime$.
\begin{Prop}\label{QA}
A sub-pseudogroup $\cH$ of $\bG(P)$ is quasi-analytic if there exists an open subset $\underline{\cT}$ of $\cT$ meeting all $\cH$-orbits such that $\cH\vert_{\underline{\cT}}$ is quasi-analytic. 
\end{Prop}
\Preuve For every $\psi\!\in\!\cH$, suppose that $\psi$ is identity on a connected open subset $V\!\subset\! \operatorname{dom}\psi$ (Recall that in this section, we assume that $\psi\!\in\!\bG(P)$ has connected domains). Let us prove that for every $z\!\in\! \overline{V}$, $\psi$ is the identity on an open neighborhood of $z$. 
Since $\underline{\cT}$ meets all $\cH$-orbits, there is $h\!\in\! \cH$ such that $z\!\in\!\operatorname{dom}h$ and $h(z)\!\in\! \underline{\cT}$.\smallskip\\
As $\underline{\cT}$ is open, we take $\epsilon\!>\!0$ such that the ball $B\big(h(z),\epsilon\big)$ at $h(z)$ with radius $\epsilon$ is contained in $\underline{\cT}$. By continuity of $h$ at $z$, there exists $\delta\!>\!0$ such that for every $w\!\in\!B(z,\delta)$, $d\big(h(z),h(w)\big)\!<\!\epsilon$. Hence, $h(w)\!\in\!B\big(h(z),\epsilon\big)\!\subset\! \underline{\cT}$, \emph{i.e.} $\operatorname{im}h\subset \underline{\cT}$.\smallskip\\
Consider $\widetilde{\psi}\!=\!h\psi h^{-1}\!\in \!\cH$ and notice that $\operatorname{dom}\widetilde{\psi}\!\subset\!\operatorname{dom}h^{-1}\!=\!\operatorname{im}h\!\subset\!\underline{\cT}$. Clearly, $\widetilde{\psi}$ is the identity on an open subset of $\underline{\cT}$ whose closure contains $h(z)$ and by the quasi-analyticity of $\cH_{\underline{\cT}}$, $\widetilde{\psi}$ is the identity on an open neighborhood of $h(z)$. Therefore, $\psi$ is identity an open neighborhood of $z$.\eb
\begin{Rem}
In general, we cannot hope to apply a change of coordinates such that $P$ satisfies the condition of triangular $1$-part under new coordinates. See the following counter-example.
\end{Rem}
\begin{CExa} In case that $q\!=\!2$, let the order $1$-part of $P$ be written by $P_1\!=\!a_1\frac{\partial}{\partial y}\!+\!a_2\frac{\partial}{\partial z}$ at some point with coordinate $(y,z)$ containing $(0,0)$. Let $(\tilde{y},\tilde{z})\!=\!\phi(y,z)$ be a change of coordinates and denote by $\tilde{P}_1\!=\!b_1\frac{\partial}{\partial \tilde{y}}\!+\!b_2\frac{\partial}{\partial \tilde{z}}$ the operator $P_1$ under these new coordinates. Notice that
$$
\begin{bmatrix}b_1(\phi)\\b_2(\phi)\end{bmatrix}\!=\!J(\phi)\begin{bmatrix}a_1\\a_2\end{bmatrix},
$$ 
where $J(\phi)\!=\!\begin{bmatrix}\phi^1_y&\phi^1_z\\ \phi^2_y&\phi^2_z\end{bmatrix}$ is the Jacobian of $\phi$ with non-zero determinant everywhere. The Jacobian of $(\tilde{y},\tilde{z})\mapsto \begin{bmatrix}b_1(\tilde{y},\tilde{z})\\b_2(\tilde{y},\tilde{z})\end{bmatrix}$ is given at $\phi(y,z)$ by
$$
\begin{bmatrix}b_{11}(\phi)&b_{12}(\phi)\\ b_{21}(\phi)&b_{22}(\phi)\end{bmatrix}\!=\!J(\phi)\begin{bmatrix}a_{11}&a_{12}\\a_{21}&a_{22}\end{bmatrix} J(\phi)^{-1}+\begin{bmatrix}\phi^1_{yy}a_1\!+\!\phi^1_{zy}a_2&\phi^1_{yz}a_1+\phi^1_{zz}a_2\\ \phi^2_{yy}a_1+\phi^2_{zy}a_2&\phi^2_{zy}a_1+\phi^2_{zz}a_2\end{bmatrix}J(\phi)^{-1}.
$$
Take $a_1\!=\!-z$, $a_2\!=\!y$, then $\small\begin{bmatrix}a_{11}&a_{12}\\ a_{21}&a_{22}\end{bmatrix}\!\!=\!\!\begin{bmatrix}0&-1\\1&0\end{bmatrix}$ is non-triangular.\\
If $\small\begin{bmatrix}b_{11}(\phi)&b_{12}(\phi)\\ b_{21}(\phi)&b_{22}(\phi)\end{bmatrix}$ is triangular, by calculation, either
$$
(\phi^1_y)^2\!+\!(\phi^1_z)^2\!+\!y(\phi^1_z\phi^1_{zy}\!-\!\phi^1_y\phi^1_{zz})\!+\!z(\phi^1_y\phi^1_{yz}\!-\!\phi^1_z\phi^1_{yy})\!=\!0\ 
$$
or
$$
(\phi^2_y)^2\!+\!(\phi^2_z)^2\!+\!y(\phi^2_z\phi^2_{zy}\!-\!\phi^2_{zz}\phi^2_y)\!+\!z(\phi^2_y\phi^2_{zy}\!-\!\phi^2_z\phi^2_{yy})\!=\!0.
$$
At $(y,z)\!=\!(0,0)$, either $\phi^1_y\!=\!\phi^1_z\!=\!0$ or $\phi^2_y\!=\!\phi^2_z\!=\!0$. But this implies $|J(\phi)|\!=\!0$ at $(0,0)$.
\end{CExa}

\section{Construction of transverse metric}
A basic connection on foliated bundles is constructed by the \textbf{Average Method} in \cite{Liu}. In the present work, in order to construct a transverse metric on $(M,\cF)$, equivalently a holonomy-invariant Riemannian metric on a complete transversal of $\cF$, one applies the same idea.\\
\indent In Section 4.1, we discuss the closure of sub-pseudogroups of $\bG(P)$. Under \textbf{Condition} introduced in Section 4.2, we construct a sub-pseudogroup of $\bG(P)$ whose closure is strictly transitive on a complete transversal of $\cF$.
In Section 4.3, we complete the construction of transverse metric on $(M,\cF)$. \\
\indent Let $P\!:\!C^\infty(M)\!\rightarrow\!C^\infty(M)$ be a basic transversely elliptic differential operator of order $m\!\geq\! 2$.
\subsection{Closure of sub-pseudogroups of $\bG(P)$}\ \\ \\
\indent Let $\cT$ be a complete transversal of $\cF$ satisfying \textbf{Assumption} and $\bG(P)$ be the commuting pseudogroup. Since $\cT$ is locally connected, the quasi-analyticity is equivalent to the strongly quasi-analyticity introduced in \cite[Definition 2.18]{AG}. For the Compact-Open topology, see \cite[Chapter 2]{Hir}.\\
\indent In \cite[Theorem 12.1]{AC}, also \cite[Theorem 2.34]{AG}, the closure of a compactly generated, equicontinuous and (strongly) quasi-analytic pseudogroup on a locally compact Polish space \cite[Definition 3.1]{Kec} is well-defined. In our setting, the complete transversal $\cT$ is a locally compact Polish space and by Propositions \ref{Equicontinuity}, \ref{AtoQA}, \ref{QA}, the equicontinuity and the quasi-analyticity of a sub-pseudogroup of $\bG(P)$ are acquired. Next proposition is a combination of our previous work.
\begin{Prop}\label{DefHBar}
Let $\cH$ be an \textbf{compactly generated} sub-pseudogroup of $\bG(P)$. If $P$ is analytic with triangular $1$-part on an open subset $\underline{\cT}$ of $\cT$ meeting all $\cH$-orbits, then the closure $\overline{\cH}$ of $\cH$ is a well-defined pseudogroup. Each element $\psi\!\in\!\overline{\cH}$ is a local transformation satisfying the following property:\\
for every $z\!\in\!\operatorname{dom}\psi$, there exists an open neighborhood $\mathcal{O}_z$ of $z$ in $\operatorname{dom}h$ such that the restriction $\psi\vert_{\mathcal{O}_z}$ is in the closure of $C^\infty(\mathcal{O}_z,\cT)\cap \cH$ in $C^\infty_{\mathrm{c}\text{-}\mathrm{o}}(\mathcal{O}_z,\cT)$, where $C^\infty_{\mathrm{c}\text{-}\mathrm{o}}$ representing space of smooth maps equipped with the Compact-Open topology. 
\end{Prop}
\Preuve By Proposition \ref{Equicontinuity}, \ref{AtoQA}, \ref{QA}, $\cH$ is equicontinuous and quasi-analytic. Noticed that in our setting all maps are $C^\infty$ and $\overline{\cH}$ is well-defined with respect to the Compact-Open topology of $C^\infty$ maps. We can then conclude by applying \cite[Theorem 2.34]{AG}. \eb

\begin{Prop}
Under the condition of Proposition \ref{DefHBar}, $\overline{\cH}$ is a sub-pseudogroup of $\bG(P)$. Furthermore, $\overline{\cH}$ is quasi-analytic.
\end{Prop}
\Preuve At each point $z\!\in\!\operatorname{dom}\psi$, by definition, there exists an open neighborhood $\mathcal{O}_z\!\subset\!\operatorname{dom}\psi$ such that $\psi\vert_{\mathcal{O}_z}$ is a local transformation belonging to the closure of $C^\infty_{\mathrm{c}\text{-}\mathrm{o}}(\mathcal{O}_z,\cT)\cap \cH$.\\
We may consider $\mathcal{O}_z$ and $V\!\!=\!\!\psi(\mathcal{O}_z)$ being located in local charts $(O_z,g)$, $(V,h)$. Let $y\!=\!(y_1,\cdots,y_q)$ be coordinate in $g(\mathcal{O}_z)$. Take an open neighborhood $\mathcal{O}_z^\prime$ of $z$. Clearly, for each compact set $\mathcal{K}$ with $\mathcal{O}_z^\prime\!\subset\! \mathcal{K} \subset \mathcal{O}_z$, one has $\psi(\mathcal{K})\subset V$. Let $\epsilon>0$ and recall that a subbasis of neighborhood \cite[Page 35]{Hir} is given by the set of smooth maps $\phi\!:\!\mathcal{O}_z\rightarrow V$ such that $\phi(\mathcal{K})\!\subset\! V$ and 
$$
\big|\!\big| D^k(h\phi g^{-1})\vert_y-D^k(h\psi g^{-1})\vert_y\big|\!\big|<\epsilon
$$
for all $y\!\in\! g(\mathcal{K})$, $k\!=\!0,\cdots,m$ where $m$ is the order of $P$.\\
For every $n\in\mathbb{N}^*$, take $\epsilon\!=\!\frac{1}{n}$, we obtain a sequence $\psi_n\!:\!O_z\rightarrow V$ satisfying $\psi_n(\mathcal{K})\subset V$ and for all $y\!\in\! g(\mathcal{K})$, $k\!=\!0,\cdots,m$
\begin{equation}\label{eq:6}
\big|\!\big|D^k(h\psi_ng^{-1})\vert_y-D^k(h\psi g^{-1})\vert_y\big|\!\big|<\frac{1}{n}
\end{equation}
Let us show that $\psi$ and $P$ commute on $\mathcal{O}_z^\prime$. For every $f\in C^\infty(V)$, on one hand, take $k\!=\!0$ in (\ref{eq:6}) and we get $\psi_n\rightarrow \psi$ on $\mathcal{O}_z^\prime$. Therefore, $\psi_n^*Pf\rightarrow\psi^*Pf$ is followed from
$$
(\psi_n^*Pf)\vert_w\!=\!(Pf)\vert_{\psi_n(w)}\!\rightarrow\! (Pf)\vert_{\psi(w)}\!=\!(\psi^*Pf)\vert_w.
$$
On the other hand, let the function $f\psi_n g^{-1}$ in $g(\mathcal{O}_z)$ be written by $(\alpha^1_n,\cdots,\alpha^q_n)$. For a multi-index $s\!=\!(s_1,\cdots,s_q)$, the $s$-term $a_s\partial^s\big[(\psi_n^*f)g^{-1}\big]$ of $P(\psi^*_n f)$ is composed by the derivatives of the functions $(\alpha^1_n,\cdots,\alpha^q_n)$ up to order $|s|$. As the relation (\ref{eq:6}) gives the convergence of $\psi_n$ up to order $m\!\geq\!|s|$, one get 
$$
P(\psi^*_n f)\rightarrow P(\psi^*f).
$$
Therefore, $\psi\vert_{\mathcal{O}^\prime_z}$ commutes with $P$. Consequently, $\psi\in\bG(P)$ by arbitrariness of $z$.\\
The quasi-analyticity of $\overline{\cH}\vert_{\underline{\cT}}$ is obtained by Proposition \ref{AtoQA}.  Notice that $\underline{\cT}$ meets all $\cH$-orbits, a priori all $\overline{\cH}$-orbits. We then conclude by applying Proposition \ref{QA}.\eb 
\medskip
\subsection{Transitivity of Pseudogroup}\ \\ \\
\indent A pseudogroup $\bH$ on $\bT$ is said to be \textbf{transitive} if some orbit is dense on $\bT$, even more to be \textbf{strictly transitive} if each pair of points $(z,w)\!\in\!\bT^2$, there is $h\!\in\!\bH$ such that $h(z)\!=\!w$. Intuitively, one may consider the holonomy pseudogroup and its closure for applying the Average Method as in \cite{Liu} the closure of the \'{e}tale holonomy groupoid is considered. Unfortunately, although ``the closure of holonomy pseudogroup'' is well-defined, the lack of sufficient transitivity bothers us to do local constructions in Section 4.3.2.  \smallskip\\
\indent An example is the Kronecker foliation on the $2$-torus with compact leaf diffeomorphic to $S^1$. All orbits of the holonomy pseudogroup are discrete on transversal $\cT\!=\!S^1\!\times\!\{*\}$. The action of the holonomy pseudogroup, which coincides with the one of its closure, has no chance to extend an inner-product at $z\!\in\!\cT$ to a Riemannian metric around $z$.\\
\indent Apart from this, \textbf{the transversal $\cT$ is independent from before}. 
\begin{Def}
We say that $P$ has \textbf{constant coefficients} at $z\!\in\! M$ if there exists a foliated chart at $z$ such that the coefficients $a_s$ in expression (\ref{eq:1.1}) are constants; $P$ has constant coefficients on a subset of $M$ if $P$ has at each point.
\end{Def}
\noindent Let us impose from now on the following condition.\\ \\
\noindent\fcolorbox{black}{white}{
\begin{minipage}{0.9\textwidth}
\textbf{Condition}: There exists a complete transversal $\cT$ of $\cF$ satisfying \textbf{Assumption} where $P$ has constant coefficients.
\end{minipage}
}\\  \\
\indent The rest of Section 4.2 is aimed to define a sub-pseudogroup $\cH$ of $\bG(P)$ whose closure $\overline{\cH}$ is strictly transitive on $\cT$. By \textbf{Condition} and \textbf{Assumption} (3), each $\cT_i^\prime$ locates in a foliated chart where $P$ has constant coefficients. 
At each $z\!\in\!\cT^\prime$, in a local chart $U\!\simeq\! (x,y)$ containing $z$ where $P$ has constant coefficients, we consider the vertical subspace $\mathrm{O}$ passing through $z$. We shall denote again by $P$ the restriction of $P$ to $\mathrm{O}$. Under a linear transformation to coordinate $y$, $P$ still has constant coefficients. So, we may assume $\mathrm{O}$ covering $(-1,1)^q$ with $z\!\simeq\! 0$, see Figure 3 below.\smallskip\\
\indent Let $(a,b)\in\mathbb{R}^2$ be small and rationally independent. For $k\!\!=\!\!1,\!\cdots\!,q$, the translations 
$$
\textstyle\mathrm{T}^k_{\pm a}:(-\frac{1}{3},\frac{1}{3})^q\rightarrow (-1,1)^q,\
y\mapsto (y_1,\cdots,y_k\!\pm\!a,\cdots,y_q)
$$
and similarly defined $\mathrm{T}^k_{\pm b}$ commute with $P$ because $P$ has constant coefficients. \\
Let $\mathrm{S}_z\!=\!\big\{\mathrm{T}^k_{\pm a},\mathrm{T}^k_{\pm b}\!\mid\! k\!\!=\!\!1,\cdots,q\big\}$. The pseudogroup $\big<\mathrm{S}_z\big>$ generated by $\mathrm{S}_z$ has a dense orbit $\big<\mathrm{S}_z\big>\!\cdot\! z$ on $\mathrm{O}^\prime\!\simeq\! (-\frac{1}{3},\frac{1}{3})^q$. By holonomy, $\mathrm{O}^\prime$ (resp. $\mathrm{S}_z$) is diffeomorphic to a transversal $\mathcal{O}_z^\prime\!\subset\!\cT$ $\big($resp. $\mathcal{S}_z\!=\!\{\mathbf{T}^k_{\pm a},\mathbf{T}^k_{\pm b}\mid k\!=\!1,\cdots,q\}\!\subset\! \bG(P)$$\big)$. Therefore, the orbit $\big<\mathcal{S}_z\big>\!\cdot\! z$ is dense on $\mathcal{O}_z^\prime$.\smallskip\\
\indent The cover $\{\mathcal{O}^\prime_z\mid z\in\cT^\prime\}$ of $\cT^\prime$ admits a finite cover $\{\mathcal{O}^\prime_{z_1},\cdots,\mathcal{O}^\prime_{z_\alpha}\}$ of $\cT$ because $\cT$ is relatively compact. Let $\mathcal{S}_0$ be the set of Haefliger's cocycles $h_{ij}$ on $\cT$ as in Proposition \ref{HolFinite}. 
\begin{center}
\begin{tikzpicture}[domain=0:4]
\draw plot[smooth, tension=1] coordinates {(1.2,0.4) (1.5,1.15) (2,1.75) (2.6,2.3) (2.7,2.85)};
\draw plot[smooth] coordinates {(2.7,2.85) (4.7,2.75) (6.2,2.85)};
\draw plot[smooth, tension=1] coordinates {(4.7,0.4) (5,1.15) (5.5,1.75) (6.1,2.3) (6.2,2.85)};
\draw plot[smooth] coordinates {(1.2,0.4) (3.2,0.5) (4.7,0.4)};
\draw plot[smooth] coordinates {(1.28,0.6) (3.18,0.7) (4.72,0.6)};
\draw plot[smooth] coordinates {(1.3,0.8) (3.2,0.9) (4.8,0.8)};
\draw plot[smooth] coordinates {(1.44,1) (3.34,1.1) (4.92,1)};
\draw plot[smooth] coordinates {(1.5,1.2) (3.4,1.3) (5.05,1.2)};
\draw plot[smooth] coordinates {(1.64,1.4) (3.54,1.5) (5.18,1.4)};
\draw plot[smooth] coordinates {(1.84,1.6) (3.74,1.7) (5.32,1.6)};
\draw plot[smooth] coordinates {(2.04,1.8) (3.94,1.85) (5.52,1.8)};
\draw plot[smooth] coordinates {(2.33,2) (4.22,2.05) (5.81,2)};
\draw plot[smooth] coordinates {(2.58,2.25) (4.52,2.20) (6.05,2.25)};
\draw plot[smooth] coordinates {(2.7,2.45) (4.7,2.35) (6.2,2.45)};
\draw plot[smooth] coordinates {(2.71,2.65) (4.71,2.55) (6.21,2.65)};
\draw plot[smooth, tension=1] coordinates {(2.3,1.1) (2.5,1.6) (3,2) (3.2,2.3)};
\draw plot[smooth, tension=1] coordinates {(2.3,1.1) (3.5,1.15)};
\draw plot[smooth, tension=1] coordinates {(3.5,1.15) (3.7,1.6) (4.2,2) (4.4,2.3)};
\draw plot[smooth, tension=1] coordinates {(3.2,2.3) (4.4,2.3)};
\draw plot[smooth, tension=1] coordinates {(2.5,0.75) (2.7,1.1) (3,1.6) (3.4,1.9)};
\draw plot[smooth, tension=1] coordinates {(2.5,0.75) (3.7,0.75)};
\draw plot[smooth, tension=1] coordinates {(3.7,0.75) (3.9,1.1) (4.2,1.6) (4.6,1.9)};
\draw plot[smooth, tension=1] coordinates {(3.4,1.9) (4.6,1.9)};
\draw (1.2,-0.3) -- (1.2,-1.9) -- (2.9,-1.9) -- (2.9,-0.3) -- (1.2,-0.3);
\draw (1.2,-0.5) -- (2.9,-0.5);
\draw (1.2,-0.7) -- (2.9,-0.7);
\draw (1.2,-0.9) -- (2.9,-0.9);
\draw (1.2,-1.1) -- (2.9,-1.1);
\draw (1.2,-1.3) -- (2.9,-1.3);
\draw (1.2,-1.5) -- (2.9,-1.5);
\draw (1.2,-1.7) -- (2.9,-1.7);
\draw[->] (2.4,1.2) to [out=240,in=120] (1.7,-0.5);
\draw[thick] (2,-0.3) -- (2,-1.9);
\node[above] at (2.1,-0.3) {\tiny $[-1\!,\!1]^q$};
\draw[very thick] (2.01,-0.83) -- (2.01,-1.36);
\draw[very thick] (1.99,-0.83) -- (1.99,-1.36);
\draw[->] (3.2,-0.9) -- (2,-1.1);
\draw[thick] plot[smooth, tension=1] coordinates {(3.2,1.1) (3.4,1.6) (3.9,2) (4.1,2.3)};
\draw[thick] plot[smooth, tension=1] coordinates {(2.7,0.75) (2.9,1) (3.4,1.5) (4.1,1.9)};

\node[above] at (4,2.3) {\footnotesize$\mathrm{O}$};
\node[below] at (2.7,0.75) {\footnotesize$\mathcal{T}$};
\node[right] at (2.7,-0.9) {\tiny $[-\frac{1}{3}\!,\!\frac{1}{3}]^q$}; 
\node[right] at (3.2,1.3) {$z$};
\node at (3.6,-2) {\footnotesize Figure 3};
\end{tikzpicture}
\hspace{2em}
\begin{tikzpicture}[domain=0:4,scale=0.7]
\draw[dashed,thick] (-3,-3) -- (3,-3) -- (3,3) -- (-3,3) -- (-3,-3);
\draw[densely dotted, thick] (-1,-1) -- (1,-1) -- (1,1) -- (-1,1) -- (-1,-1);
\draw[densely dotted, thick] (-0.2,-1.6) -- (1.8,-1.6) -- (1.8,0.4) -- (-0.2,0.4) -- (-0.2,-1.6);
\draw[->] (-3.5,0) -- (3.5,0);
\draw[->] (0,-3.5) -- (0,3.5);
\draw (1,2) -- (1,1.2);
\draw (1.8,2) -- (1.8,1.2);
\draw (2,1) -- (2.8,1);
\draw (2,0.4) -- (2.8,0.4);
\node[below] at (-3,0) {$-1$};
\node[below] at (-1,0) {$-\frac{1}{3}$};
\node[below] at (3.2,0) {$1$};
\node[below] at (1.2,0) {$\frac{1}{3}$};
\node at (-2.8,2.8) {$\mathrm{O}$};
\node at (-0.8,0.8) {$\mathrm{O}^\prime$};
\node at (1.4,1.6) {$\stackrel{T^1_a}{\longrightarrow}$};
\node at (2.4,0.7) {$\downarrow$ \footnotesize $T^2_b$ \normalsize};
\node at (0,-3.7) {\footnotesize Figure 4};
\end{tikzpicture}
\end{center}
\begin{Def}\label{DefH}
Define $\cH\!=\!\big<\mathcal{S}_0\cup\mathcal{S}_{z_1}\cup\cdots\cup\mathcal{S}_{z_\alpha}\big>\vert_\cT$.
\end{Def}
\noindent Evidently, $\cH$ contains the holonomy pseudogroup on $\cT$.
\begin{Prop}\label{StrictTransitivity}
The closure $\overline{\cH}$ is well-defined and strictly transitive on $\cT$.
\end{Prop}
\Preuve Clearly, $\cH$ is compactly generated. The transitivity on $\cT$ follows from that on each $\mathcal{O}^\prime_{z_i}$. By Proposition \ref{Equicontinuity}, $\cH$ is equicontinous. As $P$ has constant coefficients, the condition in Definition \ref{Triangular1} is satisfied with zero matrix. By Proposition \ref{AtoQA}, $\cH$ is quasi-analytic and by Proposition \ref{DefHBar}, $\overline{\cH}$ is well-defined.\\ 
By \cite[Theorem 2.33]{AG}, $\cH$ is minimal, that is, all $\cH$-orbits are dense. Then for each $z\in\cT$ the orbit $\overline{\cH}\!\cdot\! z\!=\!\cT$. Therefore, $\overline{\cH}$ is strictly transitive on $\cT$.\eb

\subsection{Construction of transverse metric}\ \\ \\
Take a point $z\in\cT$. The construction contains three steps. \\[-1em]
\begin{itemize}
\item[Step 1:] Construction of $\cH$-invariant inner product on $T_z\cT$;
\item[Step 2:] Construction of $\cH$-invariant Riemannian metric around $z$;
\item[Step 3:] Construction of $\cH$-invariant Riemannian metric on $\cT$.  
\end{itemize}
\subsubsection{Construction of $\cH$-invariant inner product at $z$}\ \\
Denote the tangent space at $z$ by $V\!\!=\!\!T_z\mathcal{O}$. The space of linear maps $\mathcal{L}(V)$ is a finite dimensional normed vector space with $|\!|u|\!|\!=\!\mathrm{sup}_{|\!|v|\!|=1} |\!|u(v)|\!|$ where $|\!|v|\!|$ is induced by Riemannian metric. It is well-known that on a finite dimensional normed vector space, all norms are equivalent and induce the same topology. Equipped with the induced topology from $\mathcal{L}(V)$, the general linear group $\mathrm{GL}(V)$ is a Lie group. Let $\underline{\overline{\cH}}_z$ be the germ group of $\overline{\cH}$ at $z$ (the quotient $\underline{\overline{\cH}}$ of $\overline{\cH}$ by the germ map is a groupoid and $\underline{\overline{\cH}}_z$ is the isotropy group of $\underline{\overline{\cH}}$ at $z$, see \cite[Page 3]{Mac}, \cite[Page 5]{Mac1}, \cite[Page 112]{Moe}.) Consider the following group homomorphism
$$
j:\underline{\overline{\cH}}_z\!\rightarrow\!\mathrm{GL}(V),\ 
\underline{h}\rightarrow d_z\underline{h}, 
$$
where $d_z$ is the differential at $z$ of $\underline{h}$. Remark that $j$ is not necessarily injective.
\begin{Prop}
The closure $\overline{\Gamma}$ of $\Gamma\!\!=\!\!j(\underline{\overline{\cH}}_z)$ in $\mathcal{L}(V)$ is a compact Lie subgroup of $\mathrm{GL}(V)$.
\end{Prop}
\Preuve For every $\gamma\!\in\!\Gamma$, $|\!|\gamma|\!|\!\leq\! \mu$ by Proposition \ref{Bounds}. Thus, $\Gamma$ is a subset of $B(0,\mu)$, the ball at $0$ with radius $\mu$. Therefore, $\overline{\Gamma}$ is compact because $\overline{B(0,\mu)}$ is compact and closed. It is clear that $\overline{\Gamma}$ is a Lie subgroup of $\mathrm{GL}(V)$. \eb
\begin{Def}
Let $\big<-,-\big>$ be the inner product on $V$ and $\mu_\gamma^{\mathrm{Haar}}$ be the Haar measure on $\overline{\Gamma}$. The inner product $g_z$ on $V$ is defined by
$$
g_z\!=\!\int_{\gamma\in\overline{\Gamma}} \gamma^*\big<-,-\big> d\mu_\gamma^{\mathrm{Haar}}.
$$ 
\end{Def}
\begin{Prop}
The inner product $g_z$ is well-defined and $\underline{\overline{\cH}}_z$-invariant.
\end{Prop}
\subsubsection{Construction of Riemannian metric around $z$}\ \\
We retain notations in Section 4.2. Thanks to holonomy, it is equivalent to do the construction on $\mathrm{O}^\prime \!\simeq\! (-\frac{1}{3},\frac{1}{3})^q$. Take $\epsilon\!>\!0$ small enough. In the following, the neighborhood $\mathcal{O}\subset\cT$ is considered to be identified with $(-\frac{1}{3}\!+\!2\epsilon,\frac{1}{3}\!-\!2\epsilon)^q$. For every $y\!\in\! (-\frac{1}{3}\!+\!2\epsilon,\frac{1}{3}\!-\!2\epsilon)^q$, the translation map sending $y$ to $0$
$$
\textstyle\mathrm{T}_y\!:\!\prod_{i=1}^q (y_i\!-\!\epsilon,y_i\!+\!\epsilon)\!\rightarrow\! (-\epsilon,\epsilon)^q,\  \tilde{y}\!\mapsto\!\tilde{y}\!-\!y
$$
can be approximated by local transformations in $\big<\mathrm{S}_z\big>$. Equivalently, for $w\!\in\!\mathcal{O}$ with coordinate $y$, the local transformation $\mathbf{T}_w$ can be identified with $\mathrm{T}_y\in \overline{\cH}$.\\
\indent Let $\mathbf{I}_w:T_w\mathcal{O}\rightarrow V$
be the tangent map of $\mathbf{T}_w$. When $w$ varies on $\mathcal{O}$, $\mathbf{I}$ is regarded as a smooth section of $C^\infty(\mathcal{O},T^*\mathcal{O}\otimes \underline{V})$ where $\underline{V}\!=\!\mathcal{O}\!\times\! V$. The smoothness is inherited from the translations $\mathrm{T}_y$.
\begin{Def}
The Riemannian metric $g_{\mathcal{O}}\!\!\in\! C^\infty(\mathcal{O},S^2T^*\mathcal{O})$ is defined by $g\vert_w\!=\!\mathbf{I}_w^*g_z$.
\end{Def}
\begin{Prop}
The Riemannian metric $g_{\mathcal{O}}$ is well-defined and $\overline{\cH}\vert_{\mathcal{O}}$-invariant.
\end{Prop}
\Preuve It is well-defined because $\mathbf{I}$ is smooth. For each $h\!\in\! \overline{\cH}\vert_{\mathcal{O}}$, every $w\!\in\! \mathcal{O}$, the germ of $\widetilde{h}\!=\!\mathbf{T}_{h(w)}h\mathbf{T}_{w}^{-1}\in \overline{\cH}$ at $z$ belongs to $\underline{\overline{\cH}}_z$,
$$
h^*g\vert_{h(w)}\!=\!h^* \mathbf{I}^*_{h(w)}g_z\!=\!\mathbf{I}^*_{w} (\mathbf{I}^{-1}_{w})^*h^*\mathbf{I}^*_{h(w)}g_z\!=\!\mathbf{I}^*_{w}\big(\widetilde{h}^{-1}\big)^*g_z\!=\!\mathbf{I}_{w}^*g_z\!=\!g\vert_{w}.
$$
\vspace{-1em}\eb
\subsubsection{Global Construction of $\cH$-invariant Riemannian metric}\ \\
For $h\in\overline{\cH}$, denote by $h\mathcal{O}\!=\!h\big(\!\operatorname{dom}h\cap \mathcal{O}\big)$. The Riemannian metric $h^*g_{\mathcal{O}}$ is well-defined on the open subset $h^{-1}\mathcal{O}\subset \cT$. Notice that if $h,h^\prime\in \overline{\cH}$ with $h^{-1}\mathcal{O}\cap {h^\prime}^{-1}\mathcal{O}\neq \emptyset$, the Riemannian metrics $h^*g_{\mathcal{O}}$ and ${h^\prime}^*g_{\mathcal{O}}$ coincide on $h^{-1}\mathcal{O}\cap {h^\prime}^{-1}\mathcal{O}$ because $h^\prime h^{-1}\in \overline{\cH}\vert_{\mathcal{O}}$. By the strict transitivity of $\overline{\cH}$ on $\cT$ in Proposition \ref{StrictTransitivity}, we define the Riemannian metric $g$ on $\cT$ by gluing the Riemannian metric $h^*g_{\mathcal{O}}$ when $h$ runs over $\overline{\cH}$. Next proposition is clear.
\begin{Prop}
The Riemannian metric $g$ is well-defined and $\overline{\cH}$-invariant. In particular, $g$ is \textit{Hol}-invariant and defines a transverse metric on $(M,\cF)$, \emph{i.e.} the foliation $\cF$ is a Riemannian foliation.
\end{Prop}
\Preuve By the construction above, $g$ is clearly $\overline{\cH}$-invariant. By Definition \ref{DefH}, $\overline{\cH}$ contains the holonomy pseudogroup \textit{Hol} and thus $g$ is \textit{Hol}-invariant. By Proposition \ref{HolRiemannian}, $g$ defines a transverse metric on $(M,\cF)$, \emph{i.e.} $\cF$ is Riemannian. \eb
\noindent We summarize our work to the following theorem.
\begin{Thm}\label{Thm2}
Let $(M,\cF)$ be a compact foliated manifold and $P\!:\!C^\infty(M)\!\rightarrow\! C^\infty(M)$ be a basic transversely elliptic operator. The foliation $\cF$ is a Riemannian foliation if there exists a complete transversal $\cT$ of $\cF$ satisfying \textbf{Assumption} where $P$ has constant coefficients. 
\end{Thm}
\begin{Rem}
In \cite[Question 2.9.2]{Ka}, A.El Kacimi-Alaoui assumed that the foliation $\cF$ on $M$ admits non-constant basic functions. In fact, if $\cF$ admits only constant basic functions, equivalently $\cF$ has a dense leaf, the holonomy pseudogroup is minimal and all leaves are dense on $M$. In this case, we have the following corollary.
\end{Rem}
\begin{Cor} 
If $\cF$ has a dense leaf and $P$ has constant coefficient at some point of $M$, then $\cF$ is a Riemannian foliation.
\end{Cor}
\Preuve If $P$ has constant coefficients in the local chart $U\!\simeq\!(x,y)$ at $z\!\in\! M$, then take $\cT^\prime$ to be the vertical subspace passing through $z$ and take $\cT$ to be a transversal at $z$ located in $\cT^\prime$. Clearly, $\cT$ satisfies \textbf{Assumption}. Since $\cF$ is minimal, $\cT$ can be taken arbitrarily small and is complete. \eb
Let us give an example where the condition of Theorem \ref{Thm2} is satisfied.
\begin{Exa}
Let $(M,\cF)$ be a compact foliation. If $M$ admits a finite cover of foliated chart $\{U_i\}_{i\in I}$ such that $P$ has constant coefficients on each $U_i$, then followed the construction of $\cT$ in Proposition \ref{ConstructionTransversal} it is clear that $\cT$ satisfies \textbf{Assumption} where $P$ has constant coefficients. Therefore, $\cF$ is a Riemannian foliation.
\end{Exa}
Let us see the case of foliation constructed by suspension \cite[Page 76]{CC}. This example shows that the present work is a \textbf{generalization} of \cite{Fu}. We assume that all maps are smooth.
\begin{Exa} 
Let $B$ be a compact connected smooth manifold and $F$ be a \textbf{compact} smooth manifold $($both are without boundary$)$. Let $\widetilde{B}$ be the universal covering of $B$ and the foliation $\cF$ on $M\!=\!\Gamma\backslash (\widetilde{B}\times F)$ is given by a group homomorphism $h:\Gamma\!\rightarrow\! \mathrm{Diff}^\infty(F)$ where $\mathrm{Diff}^\infty(F)$ is the space of smooth diffeomorphisms of $F$ equipped with the Compact-Open topology. Notice that in case that $F$ is compact, the Compact-Open topology and the $C^\infty$-topology coincide. It is well-known that $\Gamma$ can be identified to $\pi_1(B,x_0)$ when a base point $x_0$ is chosen.\smallskip\\ 
\indent Clearly, $\cT\!=\!\Gamma\backslash(\{*\}\!\times\! F)$ is a transversal of $\cF$ satisfying \textbf{Assumption} and the basic transversely elliptic operator $P\!:\!C^\infty(M)\!\rightarrow\! C^\infty(M)$ is restricted to an elliptic operator $P\!:\!C^\infty(\cT)\!\rightarrow\! C^\infty(\cT)$.\smallskip\\
\indent We know that, in the suspension case, the holonomy of $\cF$ is completely characterized by the group of transformation $h(\Gamma)\!\subset\!\mathrm{Diff}^\infty(F)$. Notice that each element of $h(\Gamma)$ is a \textbf{global} transformation of $F$. Therefore, neither local constructions as in Sections 4.3.2 are needed nor local conditions of $P$ as the ``triangular 1-part'' in Definition \ref{Triangular1} and the ``constant coefficients'' in \textbf{\textrm{Condition}} are necessary. \smallskip\\
\indent By smoothness, the closure $\overline{h(\Gamma)}$ of ${h(\Gamma)}$ in $\mathrm{Diff}^\infty(F)$ also commutes with $P$. By the ellipticity of $P$ and by Furutani's result in \cite[Theorem 2.1]{Fu}, with the Compact-Open topology, $\overline{h(\Gamma)}$ is a compact Lie transformation group of $F$. A $\overline{h(\Gamma)}$-invariant $($a priori $h(\Gamma)$-invariant$)$ Riemannian metric on $F$ can be constructed by averaging an arbitrary one. Therefore, $\cF$ is a Riemannian foliation. 
\end{Exa}
\end{spacing}

\end{document}